\documentclass{amsart}

\usepackage{amsfonts,amssymb,verbatim,amsmath,amsthm,latexsym,textcomp,amscd}
\usepackage{latexsym,amsfonts,amssymb,epsfig,verbatim}
\usepackage{amsmath,amsthm,amssymb,latexsym,graphics,textcomp}
\usepackage{graphicx}
\usepackage{color}
\usepackage{url}

\input xy
\xyoption{all}

\begin{document}

\newtheorem{theorem}{Theorem}[section]
\newtheorem{prop}[theorem]{Proposition}
\newtheorem{lemma}[theorem]{Lemma}
\newtheorem{cor}[theorem]{Corollary}
\newtheorem{definition}[theorem]{Definition}
\newtheorem{conj}[theorem]{Conjecture}
\newtheorem{rmk}[theorem]{Remark}
\newtheorem{claim}[theorem]{Claim}
\newtheorem{defth}[theorem]{Definition-Theorem}

\newcommand{\boundary}{\partial}
\newcommand{\C}{{\mathbb C}}
\newcommand{\integers}{{\mathbb Z}}
\newcommand{\natls}{{\mathbb N}}
\newcommand{\ratls}{{\mathbb Q}}
\newcommand{\bbR}{{\mathbb R}}
\newcommand{\proj}{{\mathbb P}}
\newcommand{\lhp}{{\mathbb L}}
\newcommand{\tube}{{\mathbb T}}
\newcommand{\cusp}{{\mathbb P}}
\newcommand\AAA{{\mathcal A}}
\newcommand\BB{{\mathcal B}}
\newcommand\CC{{\mathcal C}}
\newcommand\DD{{\mathcal D}}
\newcommand\EE{{\mathcal E}}
\newcommand\FF{{\mathcal F}}
\newcommand\GG{{\mathcal G}}
\newcommand\HH{{\mathcal H}}
\newcommand\II{{\mathcal I}}
\newcommand\JJ{{\mathcal J}}
\newcommand\KK{{\mathcal K}}
\newcommand\LL{{\mathcal L}}
\newcommand\MM{{\mathcal M}}
\newcommand\NN{{\mathcal N}}
\newcommand\OO{{\mathcal O}}
\newcommand\PP{{\mathcal P}}
\newcommand\QQ{{\mathcal Q}}
\newcommand\RR{{\mathcal R}}
\newcommand\SSS{{\mathcal S}}
\newcommand\TT{{\mathcal T}}
\newcommand\UU{{\mathcal U}}
\newcommand\VV{{\mathcal V}}
\newcommand\WW{{\mathcal W}}
\newcommand\XX{{\mathcal X}}
\newcommand\YY{{\mathcal Y}}
\newcommand\ZZ{{\mathcal Z}}
\newcommand\CH{{\CC\HH}}
\newcommand\PEY{{\PP\EE\YY}}
\newcommand\MF{{\MM\FF}}
\newcommand\RCT{{{\mathcal R}_{CT}}}
\newcommand\PMF{{\PP\kern-2pt\MM\FF}}
\newcommand\FL{{\FF\LL}}
\newcommand\PML{{\PP\kern-2pt\MM\LL}}
\newcommand\GL{{\GG\LL}}
\newcommand\Pol{{\mathcal P}}
\newcommand\half{{\textstyle{\frac12}}}
\newcommand\Half{{\frac12}}
\newcommand\Mod{\operatorname{Mod}}
\newcommand\Area{\operatorname{Area}}
\newcommand\ep{\epsilon}
\newcommand\hhat{\widehat}
\newcommand\Proj{{\mathbf P}}
\newcommand\U{{\mathbf U}}
 \newcommand\Hyp{{\mathbf H}}
\newcommand\D{{\mathbf D}}
\newcommand\Z{{\mathbb Z}}
\newcommand\R{{\mathbb R}}
\newcommand\Q{{\mathbb Q}}
\newcommand\E{{\mathbb E}}
\newcommand\til{\widetilde}
\newcommand\length{\operatorname{length}}
\newcommand\tr{\operatorname{tr}}
\newcommand\gesim{\succ}
\newcommand\lesim{\prec}
\newcommand\simle{\lesim}
\newcommand\simge{\gesim}
\newcommand{\simmult}{\asymp}
\newcommand{\simadd}{\mathrel{\overset{\text{\tiny $+$}}{\sim}}}
\newcommand{\ssm}{\setminus}
\newcommand{\diam}{\operatorname{diam}}
\newcommand{\pair}[1]{\langle #1\rangle}
\newcommand{\T}{{\mathbf T}}
\newcommand{\inj}{\operatorname{inj}}
\newcommand{\pleat}{\operatorname{\mathbf{pleat}}}
\newcommand{\short}{\operatorname{\mathbf{short}}}
\newcommand{\vertices}{\operatorname{vert}}
\newcommand{\collar}{\operatorname{\mathbf{collar}}}
\newcommand{\bcollar}{\operatorname{\overline{\mathbf{collar}}}}
\newcommand{\I}{{\mathbf I}}
\newcommand{\tprec}{\prec_t}
\newcommand{\fprec}{\prec_f}
\newcommand{\bprec}{\prec_b}
\newcommand{\pprec}{\prec_p}
\newcommand{\ppreceq}{\preceq_p}
\newcommand{\sprec}{\prec_s}
\newcommand{\cpreceq}{\preceq_c}
\newcommand{\cprec}{\prec_c}
\newcommand{\topprec}{\prec_{\rm top}}
\newcommand{\Topprec}{\prec_{\rm TOP}}
\newcommand{\fsub}{\mathrel{\scriptstyle\searrow}}
\newcommand{\bsub}{\mathrel{\scriptstyle\swarrow}}
\newcommand{\fsubd}{\mathrel{{\scriptstyle\searrow}\kern-1ex^d\kern0.5ex}}
\newcommand{\bsubd}{\mathrel{{\scriptstyle\swarrow}\kern-1.6ex^d\kern0.8ex}}
\newcommand{\fsubeq}{\mathrel{\raise-.7ex\hbox{$\overset{\searrow}{=}$}}}
\newcommand{\bsubeq}{\mathrel{\raise-.7ex\hbox{$\overset{\swarrow}{=}$}}}
\newcommand{\tw}{\operatorname{tw}}
\newcommand{\base}{\operatorname{base}}
\newcommand{\trans}{\operatorname{trans}}
\newcommand{\rest}{|_}
\newcommand{\bbar}{\overline}
\newcommand{\UML}{\operatorname{\UU\MM\LL}}
\newcommand{\EL}{\mathcal{EL}}
\newcommand{\tsum}{\sideset{}{'}\sum}
\newcommand{\tsh}[1]{\left\{\kern-.9ex\left\{#1\right\}\kern-.9ex\right\}}
\newcommand{\Tsh}[2]{\tsh{#2}_{#1}}
\newcommand{\qeq}{\mathrel{\approx}}
\newcommand{\Qeq}[1]{\mathrel{\approx_{#1}}}
\newcommand{\qle}{\lesssim}
\newcommand{\Qle}[1]{\mathrel{\lesssim_{#1}}}
\newcommand{\simp}{\operatorname{simp}}
\newcommand{\vsucc}{\operatorname{succ}}
\newcommand{\vpred}{\operatorname{pred}}
\newcommand\fhalf[1]{\overrightarrow {#1}}
\newcommand\bhalf[1]{\overleftarrow {#1}}
\newcommand\sleft{_{\text{left}}}
\newcommand\sright{_{\text{right}}}
\newcommand\sbtop{_{\text{top}}}
\newcommand\sbot{_{\text{bot}}}
\newcommand\sll{_{\mathbf l}}
\newcommand\srr{_{\mathbf r}}
\newcommand\geod{\operatorname{\mathbf g}}
\newcommand\mtorus[1]{\boundary U(#1)}
\newcommand\A{\mathbf A}
\newcommand\Aleft[1]{\A\sleft(#1)}
\newcommand\Aright[1]{\A\sright(#1)}
\newcommand\Atop[1]{\A\sbtop(#1)}
\newcommand\Abot[1]{\A\sbot(#1)}
\newcommand\boundvert{{\boundary_{||}}}
\newcommand\storus[1]{U(#1)}
\newcommand\Momega{\omega_M}
\newcommand\nomega{\omega_\nu}
\newcommand\twist{\operatorname{tw}}
\newcommand\modl{M_\nu}
\newcommand\MT{{\mathbb T}}
\newcommand\Teich{{\mathcal T}}
\renewcommand{\Re}{\operatorname{Re}}
\renewcommand{\Im}{\operatorname{Im}}

\title{Hyperbolic p-barycenters, circumcenters, and Moebius maps}

\author{Kingshook Biswas}
\address{Indian Statistical Institute, Kolkata, India. Email: kingshook@isical.ac.in}

\begin{abstract} Given a Moebius homeomorphism $f : \partial X \to \partial Y$
between boundaries of proper, geodesically complete CAT(-1) spaces $X,Y$, and a family of probability measures $\{ \mu_x \}_{x \in X}$ on
$\partial X$, we describe a continuous family of extensions
$\{\hat{f}_p : X \to Y \}_{1 \leq p \leq \infty}$ of $f$, called the hyperbolic $p$-barycenter maps of $f$. If all the measures
$\mu_x$ have full support then for $p = \infty$ the map $\hat{f}_{\infty}$ coincides with the circumcenter map $\hat{f}$ defined previously
in \cite{biswas5}. We use this to show that if $X, Y$ are complete, simply connected manifolds with sectional curvatures $K$ satisfying
$-b^2 \leq K \leq -1$, then the circumcenter maps of $f$ and $f^{-1}$ are $\sqrt{b}$-bi-Lipschitz homeomorphisms which are inverses of each
other. It follows that closed negatively curved manifolds with the same marked length spectrum are bi-Lipschitz homeomorphic.
\end{abstract}

\bigskip

\maketitle

\tableofcontents

\section{Introduction}

\medskip

The Mostow Rigidity Theorem asserts that for $n \geq 3$ any isomorphism between fundamental groups of closed, hyperbolic $n$-manifolds is
induced by an isometry between the manifolds, so such manifolds are determined upto isometry by their
fundamental groups. For general closed negatively curved manifolds, one may ask to what extent these manifolds are
determined by their fundamental groups. Cheeger showed that if two closed negatively curved manifolds have
isomorphic fundamental groups then the total spaces of the two-frame bundles are homeomorphic (\cite{hypgroups}, 8.2.P), while Gromov showed
the unit tangent bundles are homeomorphic via a homeomorphism that preserves the orbits of the geodesic flows (\cite{hypgroups}, 8.3.E).
Farrell and Jones showed that in dimensions $n \geq 5$, the manifolds themselves must be homeomorphic (\cite{farrelljones1}, but
also gave examples of manifolds which are homeomorphic but not diffeomorphic (\cite{farrelljones2}).

\medskip

Burns and Katok conjectured that the data of fundamental group together with lengths of closed geodesics, namely
the marked length spectrum, should be enough to determine a closed negatively curved manifold up to isometry. We
recall that the marked length spectrum of such a manifold is the function $l_X : \pi_1(X) \to \R^+$ which assigns to
each based loop the length of the unique closed geodesic in its free homotopy class. Two manifolds $X, Y$ are said
to have the same marked length spectrum if there is an isomorphism $\phi: \pi_1(X) \to \pi_1(Y)$ such that $l_X = l_Y \circ \phi$.
Otal (\cite{otal2}) showed that in dimension two this implies that the manifolds $X, Y$ are isometric. While the marked
length spectrum rigidity problem remains open in higher dimensions, Hamenstadt (\cite{hamenstadt1})
showed that equality of the marked length spectrum
is equivalent to the geodesic flows of $X, Y$ being topologically conjugate.

\medskip

These problems make sense in the more general context of group actions on CAT(-1) spaces.
Bourdon showed in \cite{bourdon1}, that for a
Gromov-hyperbolic group $\Gamma$ with two quasi-convex actions on
CAT(-1) spaces $X, Y$, the natural $\Gamma$-equivariant
homeomorphism $f$ between the limit sets $\Lambda X, \Lambda Y$ is Moebius if and only if
there is a $\Gamma$-equivariant conjugacy of the abstract geodesic
flows on $\mathcal{G}\Lambda X$ and $\mathcal{G} \Lambda Y$ compatible with
$f$, where by a Moebius map we mean a map between boundaries which preserves cross-ratios.
In particular for $\tilde{X}, \tilde{Y}$ the universal covers of two closed
negatively curved manifolds $X, Y$ (with sectional curvatures bounded above by $-1$), the geodesic flows of $X, Y$ are topologically conjugate if
and only if the induced equivariant boundary map $f : \partial \tilde{X}
\to \partial \tilde{Y}$ is Moebius, both of these conditions being equivalent to equality of the
marked length spectra of $X, Y$.

\medskip

Bourdon showed (\cite{bourdon2}) that if $X$ is a rank one symmetric
space of noncompact type with maximum of sectional curvatures equal to -1 and $Y$ a CAT(-1) space then any Moebius embedding $f :
\partial X \to \partial Y$ extends to an isometric embedding $F : X \to Y$. In \cite{biswas3} the problem
of extending Moebius maps was considered for general CAT(-1) spaces, where it was shown that any Moebius homeomorphism
$f : \partial X \to \partial Y$ between boundaries of proper, geodesically complete CAT(-1) spaces $X, Y$ extends
to a $(1, \log 2)$-quasi-isometry $F : X \to Y$. In \cite{biswas5} an extension of Moebius maps was
described which is natural with respect to composition with isometries, called circumcenter extension. The circumcenter
extension $\hat{f} : X \to Y$ of a Moebius map $f$ was shown to coincide with the $(1, \log 2)$-quasi-isometric
extension described in \cite{biswas3}, and was shown to be locally $1/2$-Holder continuous.
When $X, Y$ are complete, simply connected manifolds with sectional curvatures $K$ satisfying $-b^2 \leq K \leq -1$ for some
$b \geq 1$, it was shown in \cite{biswas5} that the circumcenter map $\hat{f}$ is a $(1, (1 - \frac{1}{b})\log 2)$-quasi-isometry.

\medskip

Our main result is the following:

\medskip

\begin{theorem} \label{mainthm} Let $X, Y$ be complete, simply connected Riemannian manifolds with sectional curvatures
satisfying $-b^2 \leq K \leq -1$ for some constant $b \geq 1$. For any Moebius homeomorphism $f : \partial X \to \partial Y$
with inverse $g : \partial Y \to \partial X$, the
circumcenter extensions $\hat{f} : X \to Y$ and $\hat{g} : Y \to X$ are $\sqrt{b}$-bi-Lipschitz homeomorphisms which
are inverses of each other.
\end{theorem}

\medskip

The proof of this theorem relies on the introduction of a continuous
family of extensions $\hat{f}_p : X \to Y, 1 \leq p \leq \infty$
of the Moebius map $f$, called {\it hyperbolic p-barycenter extensions}, such that for $p = \infty$ the map $\hat{f}_{\infty}$
coincides with the circumcenter extension $\hat{f}$. The definition of these maps relies on the notion of
{\it hyperbolic p-barycenter} of a measure $\mu$ with compact support on a CAT(-1) space $X$, which is defined
to be the unique point $x \in X$ minimizing the function
$$
y \in X \mapsto ||\cosh(d(., y))||_{L^p(\mu)}
$$
This leads to a notion of {\it asymptotic hyperbolic p-barycenter} of a measure $\nu$ with compact
support on the space of geodesics $\mathcal{G} X$, obtained as a limit of hyperbolic $p$-barycenters
of measures $\mu_t$ on $X$ obtained by pushing forward the measure $\nu$ by the geodesic flow
$\phi_t : \mathcal{G} X \to \mathcal{G} X$ and the canonical projection $\pi : \mathcal{G} X \to X$.
The hyperbolic $p$-barycenter extension $\hat{f}_p$ is then defined using asymptotic hyperbolic $p$-barycenters
and the geodesic conjugacy $\phi_f : \mathcal{G} X \to \mathcal{G} Y$ induced by $f$. When $X, Y$
are manifolds, for $1 \leq p < \infty$ the maps $\hat{f}_p$ are $C^1$, and estimates on derivatives
of the maps as $p \to \infty$ eventually lead to a proof of the above theorem.

\medskip

We remark that barycenter maps have been used previously in proving rigidity results, beginning
with the work of Besson-Courtois-Gallot (\cite{bcg1}) who used it to prove their celebrated
entropy rigidity theorem for negatively curved locally symmetric spaces.
However the hyperbolic $p$-barycenter maps constructed here
do not in general coincide with the barycenter map of Besson-Courtois-Gallot.

\medskip

As an immediate corollary we have:

\medskip

\begin{theorem} \label{rigid} Let $X, Y$ be closed negatively curved manifolds of dimension $n \geq 2$. If the marked length spectra of $X, Y$
are equal, then $X, Y$ are bi-Lipschitz homeomorphic.
\end{theorem}

\medskip

Finally we remark that in \cite{biswas4} it is proved that in certain cases Moebius maps
between boundaries of simply connected negatively curved manifolds do extend to isometries
(more precisely, local and infinitesimal rigidity results are proved for deformations of the metric
on a compact set).

\medskip

\section{Preliminaries}

\medskip

We recall in this section the definitions and facts from \cite{biswas3} and \cite{biswas5}
which we will be needing.

\medskip

\subsection{Spaces of Moebius metrics}

\medskip

Let $(Z,\rho_0)$ be a compact metric space with at least four points. For a metric $\rho$ on
$Z$ the metric cross-ratio with respect to $\rho$ of a quadruple of distinct
points $(\xi, \xi', \eta, \eta')$ of $Z$ is defined by
$$
[\xi \xi' \eta \eta']_{\rho} := \frac{\rho(\xi, \eta) \rho(\xi', \eta')}{\rho(\xi,
\eta')\rho(\xi', \eta)}
$$
A diameter one metric $\rho$ on $Z$ is {\it antipodal} if for
any $\xi \in Z$ there exists $\eta \in Z$ such that $\rho(\xi,
\eta) = 1$. We assume that $\rho_0$ is diameter one and antipodal. We say two metrics $\rho_1, \rho_2$
on $Z$ are {\it Moebius equivalent} if their metric cross-ratios agree:
$$
[\xi \xi' \eta \eta']_{\rho_1} = [\xi \xi' \eta \eta']_{\rho_2}
$$
for all $(\xi, \xi', \eta, \eta')$. The space of Moebius metrics on $Z$ is defined to be the
set $\mathcal{M}(Z, \rho_0)$ of antipodal, diameter one metrics $\rho$ on Z which are
Moebius equivalent to $\rho_0$.
We will write $\mathcal{M}(Z, \rho_0) = \mathcal{M}$. We have the following from \cite{biswas3}:

%
\medskip

\begin{theorem} \label{derivative} For any $\rho_1, \rho_2 \in \mathcal{M}$, there is a positive continuous function
$\frac{d\rho_2}{d\rho_1}$ on $Z$, called the derivative of $\rho_2$ with respect to $\rho_1$, such that
the following holds (the "Geometric Mean Value Theorem"):
$$
\rho_2(\xi, \eta)^2 = \frac{d\rho_2}{d\rho_1}(\xi) \frac{d\rho_2}{d\rho_1}(\eta) \rho_1(\xi, \eta)^2
$$
for all $\xi, \eta \in Z$.

\smallskip

Moreover for $\rho_1,
\rho_2, \rho_3 \in \mathcal{M}$ we have
$$
\frac{d\rho_3}{d\rho_1} = \frac{d\rho_3}{d\rho_2} \frac{d\rho_2}{d\rho_1}
$$
and
$$
\frac{d\rho_2}{d\rho_1} = 1/\left(\frac{d\rho_1}{d\rho_2}\right)
$$
\end{theorem}

\medskip

\begin{lemma} \label{maxmin}
$$
\max_{\xi \in Z} \frac{d\rho_2}{d\rho_1}(\xi) \cdot \min_{\xi \in Z}
\frac{d\rho_2}{d\rho_1}(\xi) = 1
$$
Moreover if $\frac{d\rho_2}{d\rho_1}$ attains its maximum 
at $\xi$ and $\rho_1(\xi, \eta) = 1$ then $\frac{d\rho_2}{d\rho_1}$ attains its minimum
at $\eta$, and $\rho_2(\xi, \eta) = 1$.
\end{lemma}

\medskip

\medskip

For $\rho_1, \rho_2 \in \mathcal{M}$, we define
$$
d_{\mathcal{M}}(\rho_1, \rho_2) := \max_{\xi \in Z}
\log \frac{d\rho_2}{d\rho_1}(\xi)
$$

\medskip

From \cite{biswas3} we have:

\medskip

\begin{lemma} \label{metric} The function $d_{\mathcal{M}}$
defines a metric on $\mathcal{M}$. The metric space $(\mathcal{M}, d_{\mathcal{M}})$ is proper.
\end{lemma}

\medskip

\bigskip

\subsection{Visual metrics on the boundary of a CAT(-1) space}

\medskip

Let $(X, d_X)$ be a proper CAT(-1) space such that $\partial X$ has at least four points.

\medskip

We recall below the definitions and some
elementary properties of visual metrics and Busemann functions;
for proofs we refer to \cite{bourdon1}:

\medskip

Let $x \in X$ be a basepoint. The {\it Gromov product} of two
points $\xi, \xi' \in \partial X$ with respect to $x$ is defined by
$$
(\xi | \xi')_x = \lim_{(a,a') \to (\xi, \xi')}
\frac{1}{2}(d(x,a) + d(x,a') - d(a,a'))
$$
where $a,a'$ are points of $X$ which converge radially towards
$\xi$ and $\xi'$ respectively. The {\it visual metric} on
$\partial X$ based at the point $x$ is defined by
$$
\rho_x(\xi, \xi') := e^{-(\xi|\xi')_x}
$$
The distance $\rho_x(\xi,\xi')$ is less than or equal to one, with
equality iff $x$ belongs to the geodesic $(\xi \xi')$.

\medskip

\begin{lemma} \label{visualantipodal} If $X$ is geodesically
complete then $\rho_x$ is a diameter one antipodal metric.
\end{lemma}

\medskip

The Busemann function $B : \partial X \times X \times X \to
\mathbb{R}$ is defined by
$$
B(x, y, \xi) := \lim_{a \to \xi} d(x,a) - d(y,a)
$$
where $a \in X$ converges radially towards $\xi$.

\medskip

\begin{lemma} \label{busemann} We have $|B(x,y,\xi)| \leq
d(x,y)$ for all $\xi \in \partial x, x,y \in X$. Moreover
$B(x,y,\xi) = d(x,y)$ iff $y$ lies on the geodesic ray $[x,
\xi)$ while $B(x,y,\xi) = -d(x,y)$ iff $x$ lies on the
geodesic ray $[y, \xi)$.
\end{lemma}

\medskip

We recall the following Lemma from \cite{bourdon1}:

\medskip

\begin{lemma} \label{visualmvt} For $x, y \in X, \xi, \eta \in
\partial X$ we have
$$
\rho_y(\xi, \eta)^2 = \rho_x(\xi, \eta)^2 e^{B(x,y,\xi)} e^{B(x,y,\eta)}
$$
\end{lemma}

\medskip

An immediate corollary of the above Lemma is the following:

\medskip

\begin{lemma} \label{visualmoebius} The visual metrics $\rho_x, x \in X$
are Moebius equivalent to each other and
$$
\frac{d\rho_y}{d\rho_x}(\xi) = e^{B(x,y,\xi)}
$$
\end{lemma}

\medskip

It follows that the metric cross-ratio $[\xi\xi'\eta\eta']_{\rho_x}$ of a quadruple
$(\xi, \xi',\eta,\eta')$ is independent of the choice of $x \in
X$. Denoting this common value by $[\xi\xi'\eta\eta']$, it is
shown in \cite{bourdon2} that the cross-ratio is given by
$$
[\xi\xi'\eta\eta'] = \lim_{(a,a',b,b') \to (\xi, \xi',\eta,\eta')} \exp(\frac{1}{2}(d(a,b)+d(a',b') -
d(a,b') - d(a',b)))
$$
where the points $a,a',b,b' \in X$ converge radially towards
$\xi,\xi',\eta,\eta' \in \partial X$.

\medskip

We assume henceforth that $X$ is a proper, geodesically complete CAT(-1) space.
We let $\mathcal{M} = \mathcal{M}(\partial X, \rho_x)$ (this space is independent of the
choice of $x \in X$). From \cite{biswas3} we have:

\medskip

\begin{lemma} The map
\begin{align*}
i_X : X & \to \mathcal{M} \\
         x & \mapsto \rho_x  \\
\end{align*}
is an isometric embedding and the image is closed in $\mathcal{M}$.
\end{lemma}

\medskip

For $k > 0$ and $y,z \in X$ distinct from $x \in X$ let $\angle^{(-k^2)} y x z \in [0, \pi]$ denote the
angle at the vertex $\overline{x}$ in a comparison triangle $\overline{x}\overline{y}\overline{z}$ in the model space $\mathbb{H}_{-k^2}$ of
constant curvature $-k^2$. From \cite{biswas5} we have:

\medskip

\begin{lemma} \label{compexist} For $\xi, \eta \in \partial X$, the limit of the comparison angles $\angle^{(-k^2)} y x z$ exists as $y,z$ converge to $\xi, \eta$ along
the geodesic rays $[x,\xi), [x, \eta)$ respectively. Denoting this limit by $\angle^{(-k^2)} \xi x \eta$, it satisfies
$$
\sin \left( \frac{\angle^{(-k^2)} \xi x \eta}{2} \right) = \rho_x(\xi, \eta)^k
$$
\end{lemma}

\medskip

\begin{lemma} \label{buseform} For $x,y \in X, \xi \in \partial X$ and $k > 0$,
the limit of the comparison angles $\angle^{(-k^2)} y x z$ exists as $z$ converges to $\xi$ along
the geodesic ray $[x,\xi)$. Denoting this limit by $\angle^{(-k^2)} y x \xi$, it satisfies
$$
e^{kB(y,x,\xi)} = \cosh(kd(x,y)) - \sinh(kd(x,y)) \cos(\angle^{(-k^2)}y x \xi)
$$
\end{lemma}

\medskip


\medskip

%
%
%

\subsection{Conformal maps, Moebius maps and geodesic conjugacies}

\medskip

We recall the definitions of conformal maps, Moebius
maps, and the abstract geodesic flow of a CAT(-1) space.

\medskip

\begin{definition} A homeomorphism between metric spaces $f :
(Z_1, \rho_1) \to (Z_2, \rho_2)$ with no isolated points is said to be {\it conformal} if
for all $\xi \in Z_1$, the limit
$$
df_{\rho_1, \rho_2}(\xi) := \lim_{\eta \to \xi} \frac{\rho_2(f(\xi),
f(\eta))}{\rho_1(\xi, \eta)}
$$
exists and is positive. The positive function $df_{\rho_1,
\rho_2}$ is called the derivative of $f$ with respect to $\rho_1, \rho_2$.
We say $f$ is {\it $C^1$ conformal} if its derivative is continuous.

\medskip

Two metrics $\rho_1, \rho_2$ inducing the same topology on a set
$Z$, such that $Z$ has no isolated points,
are said to be conformal (respectively $C^1$ conformal) if the
map $id_Z : (Z, \rho_1) \to (Z, \rho_2)$ is conformal
(respectively $C^1$ conformal). In this case we denote the
derivative of the identity map by $\frac{d\rho_2}{d\rho_1}$.
\end{definition}

\medskip

\begin{definition} A homeomorphism between metric spaces $f :
(Z_1, \rho_1) \to (Z_2, \rho_2)$ (where $Z_1$ has at least four
points) is said to be Moebius if it preserves metric cross-ratios
with respect to $\rho_1, \rho_2$. The derivative of $f$ is defined
to be the derivative $\frac{df_*\rho_2}{\rho_1}$ of the Moebius
equivalent metrics $f_* \rho_2, \rho_1$ as defined in section 2
(where $f_* \rho_2$ is the pull-back of $\rho_2$ under $f$).
\end{definition}

\medskip

Any Moebius map
between compact metric spaces with no isolated points is $C^1$ conformal,
and the two definitions of the derivative of $f$ given above
coincide. Moreover any Moebius map $f$ satisfies the geometric
mean-value theorem,
$$
\rho_2(f(\xi), f(\eta))^2 = \rho_1(\xi,\eta)^2
df_{\rho_1,\rho_2}(\xi) df_{\rho_1,\rho_2}(\xi)
$$

\medskip

\begin{definition} Let $(X, d)$ be a CAT(-1) space. The abstract geodesic flow
space of $X$ is defined to be the space of bi-infinite geodesics
in $X$,
$$
\mathcal{G}X := \{ \gamma : (-\infty,+\infty) \to X | \gamma
\hbox{ is an isometric embedding} \}
$$
endowed with the topology of uniform convergence on compact
subsets. This topology is metrizable with a distance defined by
$$
d_{\mathcal{G}X}(\gamma_1, \gamma_2):= \int_{-\infty}^{\infty}
d(\gamma_1(t), \gamma_2(t)) \frac{e^{-|t|}}{2} \ dt
$$
We define also two continuous projections
\begin{align*}
\pi : \mathcal{G}X & \to X \\
             \gamma  & \mapsto \gamma(0) \\
\end{align*}
and
\begin{align*}
p : \mathcal{G}X & \to \partial X \\
             \gamma  & \mapsto \gamma(+\infty) \\
\end{align*}

It is shown in Bourdon \cite{bourdon1} that $\pi$ is
$1$-Lipschitz, while $p$ is an open mapping.

\medskip

For $x \in X$, the unit tangent sphere $T^1_x X \subset \mathcal{G}X$ is defined to be
$$
T^1_x X := \pi^{-1}(x)
$$

\medskip

The abstract geodesic flow of $X$ is defined to be the one-parameter group
of homeomorphisms
\begin{align*}
\phi_t : \mathcal{G}X & \to \mathcal{G}X \\
         \gamma       & \mapsto \gamma_t \\
\end{align*}
for $t \in \mathbb{R}$, where $\gamma_t$ is the geodesic $s
\mapsto \gamma(s+t)$.

\medskip

The flip is defined to be the map
\begin{align*}
\mathcal{F} : \mathcal{G}X & \to \mathcal{G}X \\
              \gamma & \mapsto \overline{\gamma} \\
\end{align*}
where $\overline{\gamma}$ is the geodesic $s
\mapsto \gamma(-s)$.
\end{definition}

\medskip

We observe that for a simply connected complete Riemannian manifold $X$ with
sectional curvatures bounded above by $-1$, the map
\begin{align*}
\mathcal{G}X & \to T^1 X \\
 \gamma      & \mapsto \gamma'(0) \\
\end{align*}
is a homeomorphism conjugating the abstract geodesic flow of $X$
to the usual geodesic flow of $X$ and the flip $\mathcal{F}$ to the
usual flip on $T^1 X$.

\medskip

Let $f : \partial X \to \partial Y$ be a
conformal map between the boundaries of CAT(-1) spaces $X, Y$
equipped with visual metrics. Then $f$ induces a bijection $\phi_f :
\mathcal{G}X \to \mathcal{G}Y$ conjugating the geodesic flows, which is defined as
follows:

\medskip

Given $\gamma \in \mathcal{G}X$, let
$\gamma(-\infty) = \xi, \gamma(+\infty) = \eta, x = \gamma(0)$, then there
is a unique point $y$ on the bi-infinite geodesic $(f(\xi),f(\eta))$ such that $df_{\rho_x,
\rho_y}(\eta) = 1$. Define $\phi_f(\gamma) = \gamma^*$ where
$\gamma^*$ is the unique geodesic in $Y$ satisfying $\gamma^*(-\infty) = f(\xi),
\gamma^*(+\infty) = f(\eta), \gamma^*(0) = y$. Then $\phi_f : \mathcal{G}X \to
\mathcal{G}Y$ is a bijection conjugating the geodesic flows. From \cite{biswas3} we
have:

\medskip

\begin{prop} \label{confconj} The map $\phi_f$ is a homeomorphism if $f$ is $C^1$ conformal. If $f$ is
Moebius then $\phi_f$ is flip-equivariant.
\end{prop}

\medskip

\subsection{Circumcenter extension of Moebius maps}

\medskip

Let $X$ be a proper, geodesically complete CAT(-1) space. Recall that for any bounded subset $B$ of
 $X$, there is a unique point $x$
which minimizes the function
$$
z \mapsto \sup_{y \in B} d(z, y)
$$
The point $x$ is called the {\it circumcenter} of $B$, which we will denote by $x = c(B)$.

\medskip

Given $K \leq 0$, a function $f : X \to \bbR$ is said to be $\mathcal{F}K$-convex if it is continuous and its restriction to any
geodesic satisfies $f'' + Kf \geq 0$ in the barrier sense. This means that $f \leq g$ if $g$ coincides with
$f$ at the endpoints of a subsegment and satisfies $g'' + Kg = 0$.
We have the following from \cite{alexanderbishop}:

\medskip

\begin{prop} \label{basicfns} Let $y \in X, \xi \in \partial X$. Then:

\smallskip

\noindent (1) The function $x \mapsto \cosh(d(x,y))$ is $\mathcal{F}(-1)$-convex.

\smallskip

\noindent (2) The function $x \mapsto \exp(B(x,y,\xi))$ is $\mathcal{F}(-1)$-convex.

\end{prop}

\medskip

From \cite{biswas5} we have the following two propositions:

\medskip

\begin{prop} \label{fkmin} Let $f$ be a positive, proper, $\mathcal{F}(-1)$-convex function on $X$.
Then $f$ attains its minimum at a unique point $x \in X$.
\end{prop}

\medskip

\begin{prop} \label{minconv} Let $f_n, f$ be positive, proper, $\mathcal{F}(-1)$-convex functions on $X$
such that $f_n \to f$ uniformly on compacts. If $x_n, x$ denote the points where $f_n, f$ attain their
minima, then $x_n \to x$.
\end{prop}

Let $K$ be a compact subset of $\mathcal{G}X$ such that $p(K) \subset \partial X$ is not a singleton. Define the function
$$
u_K(z) = \sup_{\gamma \in K} \exp(B(z, \pi(\gamma), \gamma(+\infty)))
$$
It is shown in \cite{biswas5} that the function $u_K$ is a proper, positive, $\mathcal{F}(-1)$-convex function, which
hence attains its minimum at a unique $x \in X$. We call this point the asymptotic circumcenter of $K$ and denote it by $x = c_{\infty}(K)$.

\medskip

The reason for the name 'asymptotic circumcenter' is explained by the following proposition from \cite{biswas5}:

\medskip

\begin{prop} \label{asympcirc} Let $K$ be a compact subset of $\mathcal{G}X$ such that $p(K) \subset \partial X$ is not a
singleton. Define for $t > 0$ bounded subsets $A_t$ of $X$ by $A_t = \pi(\phi_t(K))$,
where $\phi_t$ denotes the geodesic flow on $\mathcal{G}X$. Then
$$
c(A_t) \to c_{\infty}(K)
$$
as $t \to +\infty$, i.e. the circumcenters of the sets $A_t$ converge to the asymptotic circumcenter of $K$.
\end{prop}

\medskip

Let $f : \partial X \to \partial Y$ be a Moebius homeomorphism between boundaries of
proper, geodesically complete CAT(-1) spaces $X, Y$, and let $\phi_f : \mathcal{G}X \to \mathcal{G}Y$
denote the associated geodesic conjugacy.

\medskip

\begin{definition} The circumcenter extension of the Moebius map $f$ is the map $\hat{f} : X \to Y$
defined by
$$
\hat{f}(x) := c_{\infty}(\phi_f(T^1_x X)) \in Y
$$
(note that $p(\phi_f(T^1_x X)) = \partial Y$ is not a singleton so the asymptotic circumcenter of
$\phi_f(T^1_x X)$ exists).
\end{definition}

\medskip

In \cite{biswas3}, a $(1, \log 2)$-quasi-isometric extension $F : X \to Y$ of the Moebius map $f$ is constructed as follows.
Since $f$ is Moebius, push-forward by $f$ of metrics on $\partial X$ to metrics on $\partial Y$ gives a
map between the spaces of Moebius metrics $f_* : \mathcal{M}(\partial X) \to \mathcal{M}(\partial Y)$, which is
easily seen to be an isometry. For each $\rho \in \mathcal{M}(\partial Y)$, we can choose a nearest point to $\rho$
in the subspace of visual metrics $i_Y(Y) \subset \mathcal{M}(\partial Y)$. This defines a nearest-point projection
$r_Y : \mathcal{M}(\partial Y) \to Y$. The extension $F$ is then defined by
$$
F = r_Y \circ f_* \circ i_X
$$

\medskip

In \cite{biswas5} it is shown that if $\rho \in \mathcal{M}(\partial Y)$ is the push-forward of a visual metric on $\partial X$,
$\rho = f_* \rho_x$ for some $x \in X$, then in fact there is a unique visual metric $\rho_y \in \mathcal{M}(\partial Y)$
nearest to $\rho$, given by $y = \hat{f}(x)$, the asymptotic circumcenter of $\phi_f(T^1_x X)$. It follows that the
extension $F$ defined above is uniquely determined and equals the circumcenter extension $\hat{f}$.

\medskip

\begin{prop} \label{nearuniq} Let $x \in X$ and let $\rho = f_* \rho_x \in \mathcal{M}(\partial Y)$. Then $y = \hat{f}(x)$ is the unique
minimizer of the function $z \in Y \mapsto d_{\mathcal{M}}(\rho, \rho_z)$. In particular, $\hat{f} = F$, so $\hat{f}$ is a
$(1, \log 2)$-quasi-isometry.
\end{prop}

\medskip

The circumcenter extension has the following naturality properties with respect to composition with isometries:

\medskip

\begin{prop} \label{natural} Let $f : \partial X \to \partial Y$ be a Moebius homeomorphism.

\smallskip

\noindent (1) If $f$ is the boundary map of an isometry $F : X \to Y$ then $\hat{f} = F$.

\smallskip

\noindent (2) If $G : X \to X, H : Y \to Y$ are isometries with boundary maps $g, h$, then
$$
\widehat{h \circ f \circ g} = H \circ \hat{f} \circ G
$$
\end{prop}

%
%

\medskip

\section{Hyperbolic $p$-barycenters in CAT(-1) spaces}

\medskip

Let $X$ be a proper, geodesically complete CAT(-1) space.

\medskip

\begin{prop} \label{powerp} Let $f$ be a positive $\mathcal{F}(-1)$-convex function on $X$. Then for any $p \geq 1$, the
function $f^p$ is $\mathcal{F}(-1)$-convex.
\end{prop}

\medskip

\noindent{\bf Proof:} We identify the restriction of $f$ to a geodesic segment $\gamma : [-a,a] \to X$ with a function
$f : [-a, a] \to \R$. Let $f(-a) = c, f(a) = d$, then the functions
$$
g(x) := \frac{1}{\sinh(2a)} \left( (\sinh a)(d+c) \cosh x + (\cosh a)(d-c) \sinh x \right)
$$
and
$$
h(x) := \frac{1}{\sinh(2a)} \left( (\sinh a)(d^p+c^p) \cosh x + (\cosh a)(d^p-c^p) \sinh x \right)
$$
satisfy $g'' - g = 0, h'' - h = 0$, and $f$ agrees with $g$ at the endpoints of $[-a,a]$, while
$f^p$ agrees with $h$ at the endpoints of $[-a, a]$. Since $f$ is $\mathcal{F}(-1)$-convex,
$f(x) \leq g(x)$ on $[-a,a]$, hence for $x \in [-a,a]$,
\begin{align*}
f(x)^p & \leq g(x)^p \\
       & = \frac{1}{\sinh^p(2a)} \left( \sinh (a + x) d + \sinh (a - x) c \right)^p \\
       & = \left(\frac{\sinh(a+x) + \sinh(a-x)}{\sinh(2a)}\right)^p \left(\frac{\sinh (a + x) d + \sinh (a - x) c}{\sinh (a + x)  + \sinh (a - x)}\right)^p \\
       & \leq \left(\frac{\sinh(a+x) + \sinh(a-x)}{\sinh(2a)}\right)^p \left(\frac{\sinh (a + x) d^p + \sinh (a - x) c^p}{\sinh (a + x)  + \sinh (a - x)}\right) \\
       & \leq \left(\frac{\sinh(a+x) + \sinh(a-x)}{\sinh(2a)}\right) \left(\frac{\sinh (a + x) d^p + \sinh (a - x) c^p}{\sinh (a + x)  + \sinh (a - x)}\right) \\
       & = \frac{1}{\sinh(2a)} \left( \sinh (a + x) d^p + \sinh (a - x) c^p \right) \\
       & = h(x) \\
\end{align*}
so $f(x)^p \leq h(x)$, hence $f^p$ is $\mathcal{F}(-1)$-convex. $\diamond$

\medskip

For $x \in X$ let $v_x : X \to \R$ denote the function
$$
v_x(y) = \cosh d(x,y).
$$
Given $1 \leq p \leq \infty$ and a probability measure $\mu$ on $X$ with compact support,
define a function $u_{\mu, p}$ on $X$ by
$$
u_{\mu, p}(x) := ||v_x||_{L^p(\mu)}
$$

\begin{prop} \label{lpconvex} Suppose the support of $\mu$ is not a singleton. Then
for $1 \leq p < \infty$, the function $u^p_{\mu, p}$ is a positive, proper,
$\mathcal{F}(-1)$-convex function, while for $p = \infty$, the function $u_{\mu, \infty}$ is a positive, proper,
$\mathcal{F}(-1)$-convex function.
\end{prop}

\medskip

\noindent{\bf Proof:} Let $1 \leq p \leq \infty$. If $u_{\mu, p}(x) = 0$ for some $x$, then $v_x = 0$ $\mu$-a.e.,
but $v_x(y) > 0$ for $y \neq x$, so we must have $\mu(\{x\}) = 1$ and $\hbox{supp}(\mu) = \{x\}$, a contradiction.
Thus $u_{\mu,p}$ is a positive function.

\medskip

If $x_n \to x$ in $X$, then clearly $v_{x_n} \to v_x$ uniformly on compacts, hence
$u_{\mu, p}(x_n) \to u_{\mu, p}(x)$ since $\mu$ has compact support, so the functions $u_{\mu, p}$ are
continuous.

\medskip

Moreover if the support of $\mu$ is contained in a ball $B(x_0, R)$, then for $x$ outside $B(x_0, R)$
we have $v_x(y) \geq \cosh(d(x, x_0) - R)$ for all $y$ in $B(x_0, R)$, hence
$$
u_{\mu, p}(x) \geq \cosh(d(x, x_0) - R) \to \infty
$$
as $d(x, x_0) \to \infty$, so $u_{\mu, p}$ is proper.

\medskip

For $p < \infty$, by the previous proposition the function $x \mapsto \cosh^p d(x,y)$ is $\mathcal{F}(-1)$-convex for each $y \in Y$,
from which
it follows easily that the function
$$
u^p_{\mu, p}(x) = \int_X \cosh^p d(x,y) d\mu(y)
$$
is $\mathcal{F}(-1)$-convex. For $p = \infty$, since $v_x$ is continuous we have
$$
u_{\mu, \infty}(x) = ||v_x||_{L^{\infty}(\mu)} = \sup_{y \in \hbox{supp}(\mu)} \cosh d(x,y)
$$
so $u_{\mu, \infty}$ is $\mathcal{F}(-1)$-convex since a supremum of functions satisfying the $\mathcal{F}(-1)$-convexity
inequality also satisfies the $\mathcal{F}(-1)$-convexity inequality. $\diamond$

\medskip

As a consequence of the previous proposition, we can make the following definition:

\medskip

\begin{definition} \label{hypbary} For a probability measure $\mu$ on $X$ with compact support, the
hyperbolic $p$-barycenter of $\mu$ is defined to be the unique minimizer of the function
$$
x \mapsto \int_X \cosh^p d(x, y) d\mu(y)
$$
if $p < \infty$,  while for $p = \infty$ it is defined to be the unique minimizer of the function
$$
x \mapsto \sup_{y \in \hbox{supp}(\mu)} \cosh d(x,y)
$$
Note that the existence of a unique minimizer is obvious if supp$(\mu)$ is a singleton, while
if supp$(\mu)$ is not a singleton the existence is guaranteed by Propositions \ref{fkmin} and
\ref{lpconvex}. We denote the hyperbolic $p$-barycenter by $c^p(\mu)$.
\end{definition}

\medskip

Note that for $p = \infty$ it is easy to see that the hyperbolic $\infty$-barycenter equals the circumcenter
of the support of $\mu$,
$$
c^{\infty}(\mu) = c(\hbox{supp}(\mu))
$$

For $x \in X$, define a positive, continuous function $w_x$ on $\mathcal{G}X$ by
$$
w_x(\gamma) := \exp(B(x, \pi(\gamma), \gamma(+\infty)))
$$

Now let $\nu$ be a probability measure on $\mathcal{G}X$ with compact support $K$ such that
$p(K) \subset \partial X$ is not a singleton.
For $1 \leq p \leq \infty$, define a function $U_{\nu, p} : X \to \R$ by
$$
U_{\nu, p}(x) := ||w_x||_{L^p(\nu)}
$$

\begin{prop} \label{hypcenter} For $1 \leq p < \infty$, the function $U^p_{\nu, p}$ is a positive,
proper, $\mathcal{F}(-1)$-convex function. For $p = \infty$, the function $U_{\nu, \infty}$ is
a positive, proper, $\mathcal{F}(-1)$-convex function.
\end{prop}

\medskip

\noindent{\bf Proof:} Let $1 \leq p \leq \infty$. If $x_n \to x$ in $X$, then for $\gamma \in \mathcal{G} X$,
$$
|\log w_{x_n}(\gamma) - \log w_x(\gamma)| = |B(x_n, \pi(\gamma), \gamma(+\infty)) - B(x, \pi(\gamma), \gamma(+\infty))| \leq d(x_n, x)
$$
so $\log w_{x_n} \to \log w_x$ uniformly on $\mathcal{G} X$, hence $w_{x_n} \to w_x$ uniformly on compacts in
$\mathcal{G} X$, thus $U_{\nu, p}(x_n) \to U_{\nu, p}(x)$ since $\nu$ has compact support. Thus $U_{\nu, p}$
is continuous.

\medskip

Suppose $1 \leq p < \infty$. Then
$$
U^p_{\nu, p}(x) = \int_{\mathcal{G} X} \exp(pB(x, \pi(\gamma), \gamma(+\infty))) d\nu(\gamma)
$$
For each $\gamma \in \mathcal{G} X$, by Propositions \ref{basicfns} and \ref{powerp},
the function
$$
x \in X \mapsto \exp(pB(x, \pi(\gamma), \gamma(+\infty)))
$$
is $\mathcal{F}(-1)$-convex, so it follows from the above expression that $U^p_{\nu, p}$ is $\mathcal{F}(-1)$-convex.

\medskip

Now let $x_n$ be a sequence in $X$ tending to infinity, and suppose $U^p_{\nu,p}(x_n)$ does not
tend to $+\infty$. Passing to a subsequence, we may assume $U^p_{\nu,p}(x_n) \leq C$ for all $n$ for some
$C > 0$, and $x_n \to \xi \in \partial X$. Since $p(K) \subset \partial X$ is not a singleton,
we can choose $\eta \in p(K)$ such that $\eta \neq \xi$. Let $N \subset \partial X$ be a
compact neighbourhood of $\eta$ not containing $\xi$, and let $A = p^{-1}(N) \subset \mathcal{G} X$, so $\nu(A) > 0$
since $A$ is a neighbourhood of a point in $K = \hbox{supp}(\nu)$. Since the angles $\angle^{(-1)}yxz$
depend continuously on $x \in X$ and $y,z \in X \cup \partial X$ distinct from $x$,
there is an $\epsilon > 0$ and an $M \geq 1$ such that
$\angle^{(-1)}x_n \pi(\gamma) \gamma(+\infty) \geq \epsilon$ for all $n \geq M$ and all $\gamma \in A$ ).
Then for $n \geq M$ and $\gamma \in A \cap \hbox{supp}(\nu)$, by Lemma \ref{buseform} we have
\begin{align*}
& \exp(B(x_n, \pi(\gamma), \gamma(+\infty))) \\
& = \cosh(d(x_n, \pi(\gamma))) - \sinh(d(x_n, \pi(\gamma))) \cos(\angle^{(-1)}x_n \pi(\gamma) \gamma(+\infty)) \\
                                           & = e^{-d(x_n, \pi(\gamma))} + 2 \sinh(d(x_n, \pi(\gamma))) \sin^2 \left(\frac{\angle^{(-1)}x_n \pi(\gamma) \gamma(+\infty)}{2}\right) \\
                                           & \geq 2 \sinh(d(x_n, x_0) - R) \sin^2(\epsilon/2) \\
\end{align*}
where $x_0 \in X, R > 0$ are chosen such that $\pi(K) \subset B(x_0, R)$ and $M$
is large enough so that $d(x_n, x_0) > R$ for $n \geq M$.
It follows that for $n \geq M$
$$
U^p_{\nu, p}(x_n) \geq (2 \sinh(d(x_n, x_0) - R) \sin^2(\epsilon/2))^p \nu(A) \to +\infty
$$
This contradicts $U^p_{\nu,p}(x_n) \leq C$. Thus $U^p_{\nu,p}$ is proper.

\medskip

For $p = \infty$, since $w_x$ is continuous
$$
U_{\nu, \infty}(x) = ||w_x||_{L^{\infty}(\nu)} = \sup_{\gamma \in K} \exp(B(x, \pi(\gamma), \gamma(+\infty)))
$$
so $U_{\nu, \infty}$ coincides with the function $u_K$ defined in the previous section,
and $u_K$ is $\mathcal{F}(-1)$-convex and proper since
$p(K)$ is not a singleton. $\diamond$

\medskip

In light of the previous proposition, we can make the following definition:

\medskip

\begin{definition} Let $\nu$ be a probability measure on $\mathcal{G} X$ with compact support $K$, such that
$p(K)$ is not a singleton. For $1 \leq p \leq \infty$, the asymptotic hyperbolic $p$-barycenter of $\nu$ is
defined to be the unique minimizer of the function
$$
x \in X \mapsto ||\exp(B(x, \pi(.), p(.)))||_{L^p(\nu)}
$$
We denote it by $c^p_{\infty}(\nu)$.
\end{definition}

\medskip

Note that for $p = \infty$, the asymptotic hyperbolic $\infty$-barycenter of $\nu$ coincides
with the asymptotic circumcenter of $K = \hbox{supp}(\nu)$ (since $U_{\nu, \infty} = u_K$).
Asymptotic hyperbolic $p$-barycenters can be described as limits of hyperbolic $p$-barycenters
of measures on expanding sets:

\medskip

\begin{prop} Let $1 \leq p \leq \infty$. Let $\nu$ be a probability measure on $\mathcal{G} X$ with compact support $K$, such that
$p(K)$ is not a singleton. For $t > 0$, define probability measures $\mu_t$ on $X$ by $\mu_t = (\pi \circ \phi_t)_* \nu$
(where $\phi_t$ denotes the geodesic flow on $\mathcal{G} X$). Then
$$
c^p(\mu_t) \to c^p_{\infty}(\nu)
$$
as $t \to \infty$.
\end{prop}

\medskip

\noindent{\bf Proof:} For $p = \infty$, this follows from Proposition \ref{asympcirc}, since $\hbox{supp}{\mu_t} = \pi \circ \phi_t(K)$
(because $\nu$ has compact support and $\pi \circ \phi_t : \mathcal{G} X \to X$ is continuous and surjective).
Let $1 \leq p < \infty$. For $t > 0$, consider the proper, $\mathcal{F}(-1)$-convex functions
\begin{align*}
u^p_{\mu_t, p}(x) \cdot 2^p e^{-pt} & = \int_X \cosh^p(d(x,y)) \cdot 2^p e^{-pt} d\mu_t(y) \\
                                    & = \int_{\mathcal{G} X} \cosh^p(d(x, \pi \circ \phi_t (\gamma))) \cdot 2^p e^{-pt} d\nu(\gamma) \\
                                    & = \int_{\mathcal{G} X} \cosh^p(d(x, \gamma(t)))) \cdot 2^p e^{-pt} d\nu(\gamma) \\
\end{align*}
Since $p(K) \subset \partial X$ is not a singleton, for $t > 0$ large enough $\pi \circ \phi_t(K)$ is not a singleton
so the functions $u_{\mu_t, p}$ are positive.
Given a ball $B \subset X$, as $t \to \infty$
$$
d(x, \gamma(t)) - t \to B(x, \pi(\gamma), \gamma(+\infty))
$$
uniformly for $x \in B$ and $\gamma \in K$ (this is a standard consequence of exponential convergence
of asymptotic geodesics in CAT(-1) spaces), and hence
$$
\cosh^p(d(x, \gamma(t)))) \cdot 2^p e^{-pt} \to \exp(pB(x, \pi(\gamma), \gamma(+\infty)))
$$
uniformly for $x \in B$ and $\gamma \in K$. It follows easily that as $t \to \infty$
\begin{align*}
u^p_{\mu_t, p}(x) \cdot 2^p e^{-pt} & = \int_{\mathcal{G} X} \cosh^p(d(x, \gamma(t)))) \cdot 2^p e^{-pt} d\nu(\gamma) \\
                                    & \to \int_{\mathcal{G} X} \exp(pB(x, \pi(\gamma), \gamma(+\infty))) d\nu(\gamma) \\
                                    & = U^p_{\nu, p}(x) \\
\end{align*}
uniformly for $x \in B$. Since $u^p_{\mu_t, p} \cdot 2^p e^{-pt}$ and $U^p_{\nu, p}$ are positive, proper,
$\mathcal{F}(-1)$-convex functions with unique minimizers $c^p(\mu_t)$ and $c^p_{\infty}(\nu)$
respectively, it follows from Proposition \ref{minconv} that $c^p(\mu_t) \to c^p_{\infty}(\nu)$ as $t \to \infty$. $\diamond$

\medskip

\begin{prop} \label{pconv} Let $\nu$ be a probability measure on $\mathcal{G} X$ with compact
support $K$ such that $p(K) \subset \partial X$ is not a singleton. Then
$$
c^p_{\infty}(\nu) \to c^{\infty}_{\infty}(\nu) = c_{\infty}(K)
$$
as $p \to \infty$.
\end{prop}

\medskip

\noindent{\bf Proof:} Given a sequence $p_n$ tending to $+\infty$, let $x_n = c^{p_n}_{\infty}(\nu),
x = c^{\infty}_{\infty}(\nu)$.
Suppose the sequence $\{ x_n \}$ is unbounded. Passing to a subsequence we may assume $x_n \to \xi \in \partial X$.
As in the proof of Proposition \ref{hypcenter}, choosing $\eta \in p(K)$ distinct from $\xi$, and
$N \subset \partial X$ a compact neighbourhood of $\eta$ disjoint from $\xi$, letting $A = p^{-1}(N)$,
there are $\epsilon > 0$ and $M \geq 1$ such that for all $n \geq M$ we have
$$
U^{p_n}_{\nu, p_n}(x_n) \geq (2 \sinh(d(x_n, x_0) - R) \sin^2(\epsilon/2))^{p_n} \nu(A)
$$
where $x_0 \in X, R > 0$ are such that $B(x_0, R)$ contains $\pi(K)$. Thus
\begin{align*}
U_{\nu, p_n}(x) & \geq U_{\nu, p_n}(x_n) \\
                & \geq 2 \sinh(d(x_n, x_0) - R) \sin^2(\epsilon/2) \nu(A)^{1/p_n} \\
                & \to +\infty \\
\end{align*}
contradicting the fact that
$$
U_{\nu, p_n}(x) = ||w_x||_{L^{p_n}(\nu)} \to ||w_x||_{L^{\infty}(\nu)} = U_{\nu, \infty}(x)
$$
This proves that $\{ x_n \}$ is bounded. Let $y$ be a limit point of the sequence $\{x_n\}$.
Passing to a subsequence we may assume $x_n \to y$. Then the functions $w_{x_n}$ on $\mathcal{G} X$
converge uniformly to $w_y$ on the compact $K$, hence
$$
U_{\nu, p_n}(x_n) = ||w_{x_n}||_{L^{p_n}(\nu)} \to ||w_y||_{L^{\infty}(\nu)} = U_{\nu, \infty}(y)
$$
while
$$
U_{\nu, p_n}(x_n) \leq U_{\nu, p_n}(x) \to  U_{\nu, \infty}(x)
$$
so it follows that
$$
U_{\nu, \infty}(y) \leq U_{\nu, \infty}(x)
$$
which implies $y = x$ since $x$ is the unique minimizer of the function $U_{\nu, \infty}$. Thus the
only limit point of the bounded sequence $\{ x_n \}$ is $x$, hence $x_n \to x$ as required. $\diamond$

\medskip

\section{Hyperbolic $p$-barycenter extension of Moebius maps}

%

\medskip

\begin{definition} Let $X, Y$ be proper, geodesically complete CAT(-1) spaces, and let $f : \partial X \to \partial Y$
be a Moebius homeomorphism. Let $\phi = \phi_f : \mathcal{G} X \to \mathcal{G} Y$ be the associated
geodesic conjugacy.
Given $1 \leq p \leq \infty$, and a family of probability measures $\mathcal{M} = \{ \mu_x \}_{x \in X}$ on $\mathcal{G} X$
such that $\hbox{supp}(\mu_x) = T^1_x X$ for all $x$, the hyperbolic $p$-barycenter extension of $f$ with
respect to the family $\mathcal{M}$ is the map $\hat{f}_{p, \mathcal{M}} : X \to Y$ defined by
$$
\hat{f}_{p, \mathcal{M}}(x) := c^p_{\infty}(\phi_* \mu_x)
$$
(note that $\hbox{supp}(\phi_* \mu_x) = \phi(T^1_x X)$, and $p(\phi(T^1_x X)) = \partial Y$, so the
asymptotic hyperbolic $p$-barycenter of $\phi_* \mu_x$ exists).
\end{definition}

\medskip

Since $\hbox{supp}(\phi_* \mu_x) = \phi(T^1_x X)$, for $p = \infty$ the point $\hat{f}_{\infty, \mathcal{M}}(x)$
is the unique minimizer of the function
$$
z \in Y \mapsto \sup_{\gamma \in \phi(T^1_x X)} \exp(B(z, \pi(\gamma), p(\gamma)))
$$
hence
$$
\hat{f}_{\infty, \mathcal{M}} = \hat{f}(x)
$$
where $\hat{f} : X \to Y$ is the circumcenter extension of $f$.
It follows easily from Proposition \ref{pconv} that for all $x \in X$,
$$
\hat{f}_{p, \mathcal{M}}(x) \to \hat{f}(x)
$$
as $p \to \infty$.

\medskip

The following
proposition is straightforward, we omit the proof:

\medskip

\begin{prop} \label{equivar} If $\Gamma$ is a group acting by isometries on $X$ and $Y$, and $f : \partial X \to \partial Y$
is $\Gamma$-equivariant, and the family $\mathcal{M}$ is $\Gamma$-equivariant, i.e. $\gamma_* \mu_x = \mu_{\gamma x}$ for
all $x$, then the extension $\hat{f}_{p, \mathcal{M}} : X \to Y$ is $\Gamma$-equivariant.
\end{prop}

\medskip

For the rest of this article, we will only consider the case when $X, Y$ are complete, simply connected
manifolds with sectional curvatures $K$ satisfying $-b^2 \leq K \leq -1$. We fix two such manifolds
$X, Y$ and a Moebius homeomorphism $f : \partial X \to \partial Y$. We introduce some notation:

\medskip

We identify as usual $\mathcal{G} X$ with $T^1 X$, and the map $p : \mathcal{G} X \to \partial X$
with a map $p : T^1 X \to \partial X$. We identify the map $\phi_f : \mathcal{G} X \to \mathcal{G} Y$
with a map $\phi : T^1 X \to T^1 Y$ conjugating the geodesic flows.
For $x \in X, \xi \in \partial X$, we denote by $\overrightarrow{x\xi} \in T^1_x X$ the
tangent vector $\gamma'(0)$ where $\gamma$ is the unique unit speed geodesic such that
$\gamma(0) = x, \gamma(+\infty) = \xi$. Denote by $q_x : \partial X \to T^1_x X$ the map $\xi \mapsto \overrightarrow{x\xi}$.
Then $q_x$ is a homeomorphism with inverse given by the restriction to $T^1_x X$ of $p : T^1 X \to \partial X$.

\medskip

For $x \in X$ and $y, z \in X \cup \partial X$ distinct from $x$, we denote by $\angle yxz \in [0, \pi]$ the Riemannian
angle between the geodesic rays $[xy)$ and $[xz)$ at $x$. We note that the upper and lower bounds on
sectional curvatures imply upper and lower bounds on Riemannian angles by comparison angles,
$$
\angle^{(-b^2)} y x z \leq \angle y x z \leq \angle^{(-1)} y x z
$$

For the rest of this article, we fix a probability measure $\mu$ on $\partial X$ such that $\hbox{supp}(\mu) = \partial X$.
Define and fix a family $\mathcal{M} = \{ \mu_x, x \in X \}$ of probability measures on the unit tangent spheres $T^1_x X$ by $\mu_x := (q_x)_* \mu$,
so that $\hbox{supp}(\mu_x) = T^1_x X$. For $1 \leq p \leq \infty$, we will denote simply by $F_p : X \to Y$ the
hyperbolic p-barycenter extension $\hat{f}_{p, \mathcal{M}} : X \to Y$ of $f : \partial X \to \partial Y$ with respect to this family
$\mathcal{M}$, and we will denote by $F : X \to Y$ the circumcenter extension $\hat{f} : X \to Y$ of $f$.
We note that then by Proposition \ref{pconv},
$F_p \to F$ pointwise on $X$ as $p \to \infty$.

\medskip

Let $1 \leq p < \infty$. For any $x \in X$, the point $F_p(x) \in Y$ is the unique minimizer of the function
\begin{align*}
z \in Y & \mapsto \int_{T^1 Y} \exp(pB(z, \pi(w), p(w)) d(\phi_* \mu_x)(w) \\
                & = \int_{\partial X} \exp(pB(z, \pi(\phi(\overrightarrow{x\xi})), f(\xi))) d\mu(\xi) \\
\end{align*}
(the equality above following from $\phi_* \mu_x = (\phi \circ q_x)_* \mu$), while the point $F(x) \in Y$
is the unique minimizer of the function
$$
z \in Y \mapsto \sup_{\xi \in \partial X} \exp(B(z, \pi(\phi(\overrightarrow{x\xi})), f(\xi)))
$$
We recall some facts about Busemann functions. Given $y \in Y$ and $\eta \in \partial Y$, the Busemann function
$z \in Y \mapsto B(z, y, \eta)$ is a $C^2$ convex function on $Y$. We denote its gradient vector field by $\nabla B^{\eta}$
(it is independent of the choice of the point $y \in Y$), which is given at a point $z \in Y$
by $\nabla B^{\eta}(z) = - \overrightarrow{z\eta}$.
The Hessian of the Busemann function,
$$
w \in T_z Y \mapsto d^2 B^{\eta}_z(w,w) := < \nabla_w \nabla B^{\eta}, w >
$$
is a nonnegative definite quadratic form on $T_z Y$ which can be described in terms of unstable Jacobi fields as follows:

\medskip

Let $\gamma$ denote the unique geodesic such that $\gamma'(0) = -\overrightarrow{z\eta}$. Recall that an unstable Jacobi field is a
Jacobi field $J$ along $\gamma$ such that $||J(t)||$ is bounded for $t \leq 0$. For any $w \in T_z Y$, there exists a
unique unstable Jacobi field $J_w$ along $\gamma$ such that $J_w(0) = w$. We then have
$$
d^2 B^{\eta}_z(w,w) = < J'_w(0), J_w(0) > = \frac{1}{2} \frac{d}{dt}_{t = 0} ||J_w(t)||^2
$$
The upper and lower bounds of $-1$ and $-b^2$ on sectional curvatures give bounds on the growth of unstable Jacobi fields,
which together with the above expression gives
$$
||w||^2 \leq d^2 B^{\eta}_z(w,w) \leq b ||w||^2
$$
for $w \in \overrightarrow{z\eta}^{\bot}$, while for $w = \overrightarrow{z\eta}$ we have for any $v \in T_z Y$
$$
d^2 B^{\eta}_z(w,v) = < \nabla_w \nabla B^{\eta}, v > = 0
$$
because the integral curves of $\nabla B^{\eta}$ are geodesics (backward asymptotic to $\eta$).
It follows that for any $w \in T_z Y$, if $w^{\bot \eta}$ denotes the orthogonal projection
of $w$ to $\overrightarrow{z\eta}^{\bot}$, then
$$
||w^{\bot \eta}||^2 \leq d^2 B^{\eta}_z(w,w) \leq b ||w^{\bot \eta}||^2
$$

\medskip

\begin{prop} \label{buseconvex} Let $1 \leq p < \infty$. The function
$$
u : z \in Y \mapsto \int_{\partial X} \exp(pB(z, \pi(\phi(\overrightarrow{x\xi})), f(\xi))) d\mu(\xi)
$$
is $C^2$ and strictly convex, i.e. the Hessian $d^2 u_z$ is positive definite for all
$z \in Y$.
\end{prop}

\medskip

\noindent{\bf Proof:} For each $\xi \in \partial X$, the function $u_{\xi} : z \mapsto \exp(pB(z, \pi(\phi(\overrightarrow{x\xi})), f(\xi)))$
is $C^2$, with gradient given by $\nabla u_{\xi}(z) = p \nabla B^{f(\xi)}(z) \exp(pB(z, \pi(\phi(\overrightarrow{x\xi})), f(\xi)))$, and Hessian given
by
\begin{align*}
& (d^2 u_{\xi})_z(w,w) \\
& = (p < \nabla_w \nabla B^{f(\xi)}, w > + p^2 < \nabla B^{f(\xi)}(z), w >^2) \exp(pB(z, \pi(\phi(\overrightarrow{x\xi})), f(\xi))) \\
                       & \geq (p ||w^{\bot f(\xi)}||^2 + p^2 < \overrightarrow{z f(\xi)}, w >^2) \exp(pB(z, \pi(\phi(\overrightarrow{x\xi})), f(\xi))) \\
                       & > 0 \\
\end{align*}
for $w \neq 0$ since either $||w^{\bot f(\xi)}||^2 > 0$ or $< \overrightarrow{z f(\xi)}, w >^2 > 0$. It follows that
$u$ is $C^2$ with Hessian given by
$$
d^2 u_z(w,w) = \int_{\partial X} (d^2 u_{\xi})_z(w,w) d\mu(\xi) > 0
$$
for $w \neq 0$ since $(d^2 u_{\xi})_z(w,w) > 0$ for all $\xi \in \partial X$. $\diamond$

\medskip

\begin{prop} \label{pbalance} For $x \in X$, the point $F_p(x) \in Y$ is the unique $z \in Y$ such that
$$
\int_{\partial X} < \overrightarrow{zf(\xi)}, w > \exp(pB(z, \pi(\phi(\overrightarrow{x\xi})), f(\xi))) d\mu(\xi) = 0
$$
for all $w \in T_z Y$.
\end{prop}

\medskip

\noindent{\bf Proof:} Since $F_p(x)$ is the unique minimizer of the function $u$ of the previous
proposition, for any $w \in T_{F_p(x)} Y$
\begin{align*}
0 & = < \nabla u(F_p(x)), w > \\
  & = \int_{\partial X} p < \nabla B^{f(\xi)}(F_p(x)), w > \exp(pB(z, \pi(\phi(\overrightarrow{x\xi})), f(\xi))) d\mu(\xi) \\
  & = - p \int_{\partial X} < \overrightarrow{F_p(x)f(\xi)}, w > \exp(pB(z, \pi(\phi(\overrightarrow{x\xi})), f(\xi))) d\mu(\xi). \\
\end{align*}
Moreover since $u$ is strictly convex with unique minimizer $F_p(x)$, if $\nabla u(z) = 0$ for some $z \in Y$
then $z = F_p(x)$. $\diamond$

\medskip

\begin{definition} A probability measure $\nu$ on $\partial Y$ is said to be balanced at $z \in Y$ if
the vector-valued integral
$$
\int_{\partial Y} \overrightarrow{z\eta} d\nu(\eta) = 0 ,
$$
or equivalently if for all $w \in T_z Y$
$$
\int_{\partial Y} < \overrightarrow{z\eta}, w > d\nu(\eta) = 0 ,
$$
\end{definition}

\medskip

For $1 \leq p < \infty$ and $x \in X$, define probability measures $\mu^x_p$ on $\partial X$ by
$$
d\mu^x_p(\xi) = c^{-1}_{x,p} \cdot \exp(pB(F_p(x), \pi(\phi(\overrightarrow{x\xi})), f(\xi))) d\mu(\xi)
$$
where $c_{x,p} > 0$ is the constant defined by
$$
c_{x, p} = \int_{\partial X} \exp(pB(F_p(x), \pi(\phi(\overrightarrow{x\xi})), f(\xi))) d\mu(\xi)
$$
Then the previous proposition says that the measure $f_* \mu^x_p$ on $\partial Y$ is balanced
at the point $F_p(x) \in Y$. The following characterization of the circumcenter map
$F : X \to Y$ in terms
of balanced measures will be useful:

\medskip

\begin{prop} \label{circbal} Let $x \in X$ and $y \in Y$. Let $K \subset \partial X$ be the set
where the function
$$
\xi \in \partial X \mapsto \log \frac{df_* \rho_x}{d\rho_y}(f(\xi))
$$
attains its maximum value. Then the following are equivalent:

\medskip

\noindent (1) $y = F(x)$.

\medskip

\noindent (2) For any $w \in T_y Y$, there exists $\xi \in K$ such that $< \overrightarrow{yf(\xi)}, w > \leq 0$.

\medskip

\noindent (3) The convex hull in $T_y Y$ of the compact $\{ \overrightarrow{yf(\xi)} | \xi \in K \} \subset T^1_y Y$ contains
the origin of $T_y Y$.

\medskip

\noindent (4) There exists a probability measure $\nu$ on $\partial X$ with support contained in $K$ such that $f_* \nu$
is balanced at $y$.
\end{prop}

\medskip

\noindent{\bf Proof:} (1) $\Rightarrow$ (2): Given that $y = F(x)$,
suppose there exists $w \in T^1_y Y$ such $< \overrightarrow{yf(\xi)}, w > > 0$ for all $\xi \in K$.
Then we can choose a neighbourhood $N$ of $K$ in $\partial X$ and $\epsilon, \delta > 0$ such that $< \overrightarrow{yf(\xi)}, w > \geq \epsilon$ for
all $\xi \in N$, and
$$
\log \frac{df_* \rho_x}{d\rho_y}(f(\xi)) \leq M -\delta
$$
for $\xi \in \partial X - N$, where $M = d_{\mathcal M}(f_* \rho_x, \rho_y)$
is the maximum value of the function $\xi \in \partial X \mapsto \log \frac{df_* \rho_x}{d\rho_y}(f(\xi))$.

\medskip

Let $y_t = \exp_y(tw)$. As $t \to 0$, for $\xi \in N$ we have
\begin{align*}
\log \frac{df_* \rho_x}{d\rho_{y_t}}(f(\xi)) & = \log \frac{df_* \rho_x}{d\rho_y}(f(\xi)) - \log \frac{d\rho_{y_t}}{d\rho_y}(f(\xi)) \\
                                             & \leq M + B(y_t, y, f(\xi)) \\
                                             & = M - t < \overrightarrow{yf(\xi)}, w > + o(t) \\
                                             & \leq M - \epsilon t + o(t) \\
                                             & < M \\
\end{align*}
for $t > 0$ small enough depending only on $\epsilon$, while for $\xi \in \partial X - N$ we have
\begin{align*}
\log \frac{df_* \rho_x}{d\rho_{y_t}}(f(\xi)) & \leq (M - \delta) + B(y_t, y, f(\xi)) \\
                                             & \leq M - \delta + t \\
                                             & < M \\
\end{align*}
for $0 < t < \delta$. Thus for $t > 0$ small enough,
$$
d_{\mathcal M}(f_* \rho_x, \rho_{y_t}) = \sup_{\xi \in \partial X} \log \frac{df_* \rho_x}{d\rho_{y_t}}(f(\xi)) < M = d_{\mathcal M}(f_* \rho_x, \rho_y)
$$
contradicting the fact that $y$ is the unique minimizer of the function $p \in Y \mapsto d_{\mathcal M}(f_* \rho_x, \rho_p)$.
This proves $(1) \Rightarrow (2)$.

\medskip

\noindent $(2) \Rightarrow (1)$: Let $z \in Y$ be distinct from $y$, and let $w \in T^1_y Y$ be the initial velocity of the
geodesic joining $y$ to $z$. By hypothesis, there exists $\xi \in K$ such that $< \overrightarrow{yf(\xi)}, w > \leq 0$. Since $\overrightarrow{yf(\xi)}$
is the inward pointing normal to the boundary of the horoball $H = \{ p \in Y : B(p, y, f(\xi)) \leq 0 \}$ which is strictly
convex, it follows that $z \notin H$, so $B(z, y, f(\xi)) > 0$, hence
\begin{align*}
d_{\mathcal M}(f_* \rho_x, \rho_z) & = \sup_{\eta \in \partial X} \log \frac{df_* \rho_x}{d\rho_z}(f(\eta)) \\
                                       & \geq \log \frac{df_* \rho_x}{d\rho_z}(f(\xi)) \\
                                       & = \log \frac{df_* \rho_x}{d\rho_y}(f(\xi)) - \log \frac{d\rho_z}{d\rho_y}(f(\xi)) \\
                                       & = d_{\mathcal M}(f_* \rho_x, \rho_y) + B(z, y, f(\xi)) \\
                                       & > d_{\mathcal M}(f_* \rho_x, \rho_y) \\
\end{align*}
thus $y$ minimizes the function $p \in Y \mapsto d_{\mathcal M}(f_* \rho_x, \rho_p)$, so $y = F(x)$.

\medskip

\noindent $(2) \Rightarrow (3)$: Suppose the convex hull $L \subset T_y Y$ of the compact
$\{ \overrightarrow{yf(\xi)} | \xi \in K \} \subset T^1_y Y$ does not contain the origin of $T_y Y$, then there is a
hyperplane in $T_y Y$ separating $L$ from the origin, so if $w \in T^1_y Y$ is a unit normal vector to this
hyperplane, then $< v, w > > 0$ for all $v \in L$ (after possibly replacing $w$ by $-w$ if necessary), in
particular $< \overrightarrow{yf(\xi)}, w > > 0$ for all $\xi \in K$, a contradiction to our hypothesis.

\medskip

\noindent $(3) \Rightarrow (4)$: A convex combination $\sum_{i = 1}^k \alpha_i \overrightarrow{yf(\xi)}$ where
$\xi_1, \dots, \xi_k \in K$ can be written as a vector-valued integral $\int_{\partial X} \overrightarrow{yf(\xi)} d\lambda(\xi)$
where $\lambda$ is the probability measure supported on the finite set $\{ \xi_1, \dots, \xi_k \}$ with masses
$\alpha_1, \dots, \alpha_k$. Any point in the convex hull of the compact $\{ \overrightarrow{yf(\xi)} | \xi \in K \}$ is a limit
of such convex combinations, in particular by our hypothesis we have $\int_{\partial X} \overrightarrow{yf(\xi)} d\lambda_n(\xi) \to 0$
for some sequence of probability measures $\lambda_n$ supported on $K$, taking a weak limit of these measures gives
a probability measure $\nu$ supported on $K$ such that $\int_{\partial X} \overrightarrow{yf(\xi)} d\nu(\xi) = 0$.

\medskip

\noindent $(4) \Rightarrow (2)$: Suppose there is $w \in T_y Y$ such that $< \overrightarrow{yf(\xi)}, w > > 0$ for all $\xi \in K$,
then since $\nu$ is supported on $K$ we have $\int_{\partial X} < \overrightarrow{yf(\xi)}, w > d\nu(\xi) > 0$, a contradiction to the
fact that $f_* \nu$ is balanced at $y$.
$\diamond$

\medskip

\begin{prop} \label{implicit} The hyperbolic $p$-barycenter map $F_p : X \to Y$ is $C^1$ and its derivative satisfies
\begin{align*}
\int_{\partial X} d^2 B^{f(\xi)}_{F_p(x)}(DF_p(v), DF_p(v)) d\mu^x_p(\xi) +
p \int_{\partial X} < DF_p(v), \overrightarrow{F_p(x)f(\xi)} >^2 d\mu^x_p(\xi) \\
= p \int_{\partial X} < DF_p(v), \overrightarrow{F_p(x)f(\xi)} > < v, \overrightarrow{x\xi} > d\mu^x_p(\xi) \\
\end{align*}
for all $x \in X, v \in T_x X$.
\end{prop}

\medskip

\noindent{\bf Proof:} Note that for $x_1, x_2 \in X$ and $\xi \in \partial X$, the
geodesic conjugacy $\phi : T^1 X \to T^1 Y$ satisfies
$$
B(\pi(\phi(\overrightarrow{x_1\xi})), \pi(\phi(\overrightarrow{x_2\xi})), f(\xi)) = B(x_1, x_2, \xi)
$$
hence, fixing a basepoint $x_0 \in X$, the function $(x,y) \in X \times Y \mapsto B(y, \pi(\phi(\overrightarrow{x\xi})), f(\xi))$
can be written as
\begin{align*}
B(y, \pi(\phi(\overrightarrow{x\xi})), f(\xi)) & = B(y, \pi(\phi(\overrightarrow{x_0\xi})), f(\xi)) + B(\pi(\phi(\overrightarrow{x_0\xi})), \pi(\phi(\overrightarrow{x\xi})), f(\xi)) \\
                                    & = B(y, \pi(\phi(\overrightarrow{x_0\xi})), f(\xi)) + B(x_0, x, \xi) \\
\end{align*}
and is hence $C^2$ as a function of $(x,y)$ (even though $\phi$ is not necessarily even $C^1$).
Thus letting $e_1(y), \dots, e_n(y), y \in Y$ be a smooth orthonormal frame field on $Y$, we can define a $C^1$ function
$H = (H_1, \dots, H_n) : X  \times Y \to \R^n$ by
$$
H_i(x, y) = \int_{\partial X} < \nabla B^{f(\xi)}(y), e_i(y) > \exp(pB(y, \pi(\phi(\overrightarrow{x\xi})), f(\xi))) d\mu(\xi)
$$
Then by Proposition \ref{pbalance}, $F_p(x) \in Y$ is defined implicitly by the equation $H(x, F_p(x)) = 0$.
For $y = F_p(x)$ and $w \in T_y Y$, let $\nabla_w e_i = \sum_j \eta_{ij} e_j$, then, using $H(x,y) = 0$,
the partial derivative $(D_y H_i)(w)$ is given by
\begin{align*}
& (D_y H_i)(w) \\
& = \int_{\partial X} < \nabla_w \nabla B^{f(\xi)}, e_i(y) > \exp(pB(y, \pi(\phi(\overrightarrow{x\xi})), f(\xi))) d\mu(\xi)
+ \sum_j \eta_{ij} H_j(x,y) \\
            & + p \int_{\partial X} < \nabla B^{f(\xi)}(y), e_i(y) > < \nabla B^{f(\xi)}(y), w > \exp(pB(y, \pi(\phi(\overrightarrow{x\xi})), f(\xi))) d\mu(\xi) \\
            & = \int_{\partial X} < \nabla_w \nabla B^{f(\xi)}, e_i(y) > \exp(pB(y, \pi(\phi(\overrightarrow{x\xi})), f(\xi))) d\mu(\xi) \\
            & + p \int_{\partial X} < \nabla B^{f(\xi)}(y), e_i(y) > < \nabla B^{f(\xi)}(y), w > \exp(pB(y, \pi(\phi(\overrightarrow{x\xi})), f(\xi))) d\mu(\xi) \\
            & = \frac{1}{p} (d^2 u)_y(w, e_i(y)) \\
\end{align*}
where $u$ is the strictly convex function of Proposition \ref{buseconvex}. Since $d^2 u_y$ is positive-definite, it follows
that $D_y H$ is invertible, hence by the Implicit Function Theorem $F_p$ is $C^1$.

\medskip

Given $v \in T_x X$, let $(x_t)_{|t| < \epsilon}$ be the geodesic with initial velocity $v$, then,
as above we can write
\begin{align*}
B(F_p(x_t), \pi(\phi(\overrightarrow{x_t\xi})), f(\xi)) & = B(F_p(x_t), \pi(\phi(\overrightarrow{x\xi})), f(\xi)) + B(\pi(\phi(\overrightarrow{x\xi})), \pi(\phi(\overrightarrow{x_t\xi})), f(\xi)) \\
                                             & = B(F_p(x_t), \pi(\phi(\overrightarrow{x\xi})), f(\xi)) + B(x, x_t, \xi) \\
\end{align*}
so
$$
\frac{d}{dt}_{t = 0} B(F_p(x_t), \pi(\phi(\overrightarrow{x_t\xi})), f(\xi)) = < \overrightarrow{x\xi}, v > - < \overrightarrow{F_p(x)f(\xi)}, DF_p(v) >
$$
So differentiating the equality $H_i(x_t, F_p(x_t)) = 0$ at $t = 0$ gives, writing $\tau(\xi) = \exp(p B(F_p(x), \pi(\phi(\overrightarrow{x\xi})), f(\xi)))$,
\begin{align*}
0 & = \int_{\partial X} d^2 B^{f(\xi)}_{F_p(x)}( DF_p(v), e_i(F_p(x)) > \tau(\xi) d\mu(\xi) + \sum_j \eta_{ij} H_j(x, F_p(x)) \\
  & + p \int_{\partial X} < - \overrightarrow{F_p(x)f(\xi)}, e_i(F_p(x)) > ( < \overrightarrow{x\xi}, v > - < \overrightarrow{F_p(x)f(\xi)}, DF_p(v) > ) \tau(\xi) d\mu(\xi) \\
  & = \int_{\partial X} d^2 B^{f(\xi)}_{F_p(x)}( DF_p(v), e_i(F_p(x)) > \tau(\xi) d\mu(\xi) \\
  & + p \int_{\partial X} < - \overrightarrow{F_p(x)f(\xi)}, e_i(F_p(x)) > ( < \overrightarrow{x\xi}, v > - < \overrightarrow{F_p(x)f(\xi)}, DF_p(v) > ) \tau(\xi) d\mu(\xi) \\
\end{align*}
for $i = 1, \dots, k$. It follows that for any $w \in T_{F_p(x)}$,
\begin{align*}
0 & = \int_{\partial X} d^2 B^{f(\xi)}_{F_p(x)}( DF_p(v), w) > \tau(\xi) d\mu(\xi) \\
  & + p \int_{\partial X} < - \overrightarrow{F_p(x)f(\xi)}, w) > ( < \overrightarrow{x\xi}, v > - < \overrightarrow{F_p(x)f(\xi)}, DF_p(v) > ) \tau(\xi) d\mu(\xi) \\
\end{align*}
so the formula stated in the proposition follows by putting $w = DF_p(v)$. $\diamond$

\medskip

\begin{prop} \label{ineq1} For any $x \in X$, $v \in T_x X$,
$$
\frac{1}{p} \int_{\partial X} ||DF_p(v)^{\bot f(\xi)}||^2 d\mu^x_p(\xi) \leq \int_{\partial X} < \overrightarrow{x\xi}, v >^2 d\mu^x_p(\xi)
$$
(where $DF_p(v)^{\bot f(\xi)} \in T_{F_p(x)} Y$ denotes the orthogonal projection of $DF_p(v)$ to $\overrightarrow{F_p(x)f(\xi)}^{\bot}$).
\end{prop}

\medskip

\noindent{\bf Proof:} By the previous proposition,
\begin{align*}
& \int_{\partial X} < \overrightarrow{F_p(x)f(\xi)}, DF_p(v) >^2 d\mu^x_p(\xi) \\
& = \int_{\partial X} < \overrightarrow{F_p(x)f(\xi)}, DF_p(v) > < \overrightarrow{x\xi}, v > d\mu^x_p(\xi) \\
& - \frac{1}{p} \int_{\partial X} d^2 B^{f(\xi)}(DF_p(v), DF_p(v)) d\mu^x_p(\xi) \\
& \leq \int_{\partial X} < \overrightarrow{F_p(x)f(\xi)}, DF_p(v) > < \overrightarrow{x\xi}, v > d\mu^x_p(\xi) \\
& \leq \left( \int_{\partial X} < \overrightarrow{F_p(x)f(\xi)}, DF_p(v) >^2 d\mu^x_p(\xi) \right)^{1/2}
        \left( \int_{\partial X} < \overrightarrow{x\xi}, v >^2 d\mu^x_p(\xi) \right)^{1/2} \\
\end{align*}
hence
$$
\int_{\partial X} < \overrightarrow{F_p(x)f(\xi)}, DF_p(v) >^2 d\mu^x_p(\xi) \leq \int_{\partial X} < \overrightarrow{x\xi}, v >^2 d\mu^x_p(\xi).
$$
We then have
\begin{align*}
& \frac{1}{p} \int_{\partial X} ||DF_p(v)^{\bot f(\xi)}||^2 d\mu^x_p(\xi) \leq \frac{1}{p} \int_{\partial X} d^2 B^{f(\xi)}(DF_p(v), DF_p(v)) d\mu^x_p(\xi) \\
& = \int_{\partial X} < \overrightarrow{F_p(x)f(\xi)}, DF_p(v) > < \overrightarrow{x\xi}, v > d\mu^x_p(\xi) - \int_{\partial X} < \overrightarrow{F_p(x)f(\xi)}, DF_p(v) >^2 d\mu^x_p(\xi) \\
& \leq \left( \int_{\partial X} < \overrightarrow{F_p(x)f(\xi)}, DF_p(v) >^2 d\mu^x_p(\xi) \right)^{1/2}
        \left( \int_{\partial X} < \overrightarrow{x\xi}, v >^2 d\mu^x_p(\xi) \right)^{1/2} \\
& \leq \int_{\partial X} < \overrightarrow{x\xi}, v >^2 d\mu^x_p(\xi) \\
\end{align*}
$\diamond$

\medskip

\begin{prop} \label{maininequality} For any $x, y \in X$,
$$
\cosh(d(F_p(x), F_p(y))) \leq \int_{\partial X} \exp(B(F_p(y), F_p(x), f(\xi))) d\mu^x_p(\xi)
$$
and
$$
\cosh(bd(F_p(x), F_p(y))) \geq \int_{\partial X} \exp(bB(F_p(y), F_p(x), f(\xi))) d\mu^x_p(\xi)
$$
\end{prop}

\medskip

\noindent{\bf Proof:} Let
$$
\theta(\xi) = \angle{F_p(y)F_p(x)f(\xi)}, \theta_1(\xi) = \angle^{(-1)}{F_p(y)F_p(x)f(\xi)}, \theta_b(\xi) = \angle^{(-b^2)}{F_p(y)F_p(x)f(\xi)},
$$
then
$$
\theta_b(\xi) \leq \theta(\xi) \leq \theta_1(\xi)
$$
and by Proposition \ref{pbalance}
$$
\int_{\partial X} \cos(\theta(\xi)) d\mu^x_p(\xi) = 0
$$
since the measure $f_* \mu^x_p$ is balanced at the point $F_p(x)$. Thus
$$
\int_{\partial X} \cos(\theta_1(\xi)) d\mu^x_p(\xi) \quad \leq \quad 0 \quad \leq \quad \int_{\partial X} \cos(\theta_b(\xi)) d\mu^x_p(\xi)
$$
so the proposition follows from Lemma \ref{buseform}. $\diamond$

\medskip

We now fix a point $z \in X$. Then the following second-order Taylor expansion holds for $w \in T_{F_p(z)} Y$, uniformly in
$\xi \in \partial X$:
$$
B(\exp_{F_p(z)}(w), F_p(z), f(\xi)) = - < w, \overrightarrow{F_p(z) f(\xi)} > + \frac{1}{2}d^2 B^{f(\xi)}(w,w) + o(||w||^2)
$$
as $w \to 0$.

\medskip

We fix a unit tangent vector $v \in T^1_z X$, and consider two points $x = \exp_z(tv), y = \exp_z(-tv)$, with $t > 0$ small.

\medskip

\begin{prop} \label{estimate1} As $t \to 0$, we have
\begin{align*}
&\int_{\partial X} \exp(bB(F_p(y), F_p(x), f(\xi))) d\mu^x_p(\xi) \\
& \geq 1 + \left( 2b^2 \int_{\partial X} < DF_p(v), \overrightarrow{F_p(z)f(\xi)} >^2 d\mu^z_p(\xi) \right) t^2  \\
& + p^2 \left(\int_{\partial X} < v, \overrightarrow{z\xi} > d\mu^z_p(\xi) > \right)^2 t^2 + o(t^2) \\
\end{align*}
\end{prop}

\medskip

\noindent{\bf Proof:} Since $x = \exp_z(tv)$, we can write $F_p(x) = \exp_{F_p(z)}(t DF_p(v) + w)$ where $w = o(t)$ as $t \to 0$.
For any $\xi \in \partial X$, we have
\begin{align*}
& B(F_p(x), \pi(\phi(\overrightarrow{x\xi})), f(\xi)) \\
& = B(F_p(x), F_p(z), f(\xi)) + B(F_p(z), \pi(\phi(\overrightarrow{z\xi})), f(\xi)) +
B(\pi(\phi(\overrightarrow{z\xi})), \pi(\phi(\overrightarrow{x\xi})), f(\xi)) \\
& = B(F_p(x), F_p(z), f(\xi)) + B(F_p(z), \pi(\phi(\overrightarrow{z\xi})), f(\xi)) -
B(x, z, \xi), \\
\end{align*}
so, using the second-order Taylor expansion of Busemann functions,
\begin{align*}
& B(F_p(x), \pi(\phi(\overrightarrow{x\xi})), f(\xi)) - B(F_p(z), \pi(\phi(\overrightarrow{z\xi})), f(\xi)) \\
& = B(F_p(x), F_p(z), f(\xi)) - B(x, z, \xi) \\
& = (<v, \overrightarrow{z\xi}> - <DF_p(v), \overrightarrow{F_p(z)f(\xi)}>)t - <w, \overrightarrow{F_p(z)f(\xi)}> \\
& + \frac{1}{2}\left(d^2 B^{f(\xi)}(t DF_p(v) + w, t DF_p(v) + w) - d^2 B^{\xi}(v, v)t^2 \right) + o(t^2) \\
& = a(v, \xi)t - <w, \overrightarrow{F_p(z)f(\xi)}> + b(v, \xi)t^2 + o(t^2) \\
\end{align*}
where
\begin{align*}
a(v, \xi) & = <v, \overrightarrow{z\xi}> - <DF_p(v), \overrightarrow{F_p(z)f(\xi)}> , \\
b(v, \xi) & = \frac{1}{2}\left(d^2 B^{f(\xi)}(DF_p(v), DF_p(v) ) - d^2 B^{\xi}(v, v)\right). \\
\end{align*}
It follows that
\begin{align*}
& \frac{\exp(pB(F_p(x), \pi(\phi(\overrightarrow{x\xi})), f(\xi)))}{\exp(pB(F_p(z), \pi(\phi(\overrightarrow{z\xi})), f(\xi)))} \\
& = \exp\left(p\left( a(v, \xi)t - <w, \overrightarrow{F_p(z)f(\xi)}> + b(v, \xi)t^2 + o(t^2)\right)\right) \\
& = 1 + pa(v, \xi)t - p<w, \overrightarrow{F_p(z)f(\xi)}> + \left(pb(v, \xi) + \frac{1}{2}p^2 a(v,\xi)^2\right)t^2 + o(t^2) \\
\end{align*}
and hence
\begin{align*}
& \frac{c_{x,p}}{c_{z,p}} = \int_{\partial X} \exp(pB(F_p(x), \pi(\phi(\overrightarrow{x\xi})), f(\xi))) d\mu(\xi) \cdot c^{-1}_{z,p} \\
& = \int_{\partial X} (1 + pa(v, \xi)t - p<w, \overrightarrow{F_p(z)f(\xi)}> + \left(pb(v, \xi) + \frac{1}{2}p^2 a(v,\xi)^2\right)t^2 + o(t^2)) d\mu^z_p(\xi) \\
& = \int_{\partial X} (1 + pa(v, \xi)t + \left(pb(v, \xi) + \frac{1}{2}p^2 a(v,\xi)^2\right)t^2 + o(t^2)) d\mu^z_p(\xi) \\
\end{align*}
(where in the last line above we have used the fact that $f_* \mu^z_p$ is balanced at $F_p(z)$).

\medskip

Thus, letting $c(v,\xi) = pb(v,\xi) + \frac{1}{2}p^2 a(v,\xi)^2$, we have
\begin{align*}
& \frac{d\mu^x_p}{d\mu^z_p}(\xi) = \frac{\exp(pB(F_p(x), \pi(\phi(\overrightarrow{x\xi})), f(\xi)))}{\exp(pB(F_p(z), \pi(\phi(\overrightarrow{z\xi})), f(\xi)))} \frac{c_{z,p}}{c_{x,p}} \\
& = \frac{1 + pa(v, \xi)t - p<w, \overrightarrow{F_p(z)f(\xi)}> + \left(pb(v, \xi) + \frac{1}{2}p^2 a(v,\xi)^2\right)t^2 + o(t^2)}
{\int_{\partial X} (1 + pa(v, \xi)t + \left(pb(v, \xi) + \frac{1}{2}p^2 a(v,\xi)^2\right)t^2 + o(t^2)) d\mu^z_p(\xi)} \\
& = \left(1 + pa(v, \xi)t - p<w, \overrightarrow{F_p(z)f(\xi)}> + c(v,\xi)t^2 + o(t^2)\right) \\
& \cdot \left(1 - p\left(\int_{\partial X} a(v, \xi) d\mu^z_p(\xi)\right) \cdot t + \left(p^2 \left(\int_{\partial X} a(v, \xi)d\mu^z_p(\xi)\right)^2 -
\int_{\partial X} c(v,\xi) d\mu^z_p(\xi)\right) \cdot t^2
+ o(t^2)\right) \\
& = 1 + p\left(a(v, \xi) - \int_{\partial X} a(v, \xi) d\mu^z_p(\xi)\right)\cdot t - p<w, \overrightarrow{F_p(z)f(\xi)}> \\
& + \left(c(v,\xi) - \int_{\partial X} c(v,\xi) d\mu^z_p(\xi)
+ p^2 \left(\int_{\partial X} a(v,\xi) d\mu^z_p(\xi)\right)^2 \right)\cdot t^2 + o(t^2) \\
\end{align*}

\medskip

Letting $F_p(y) = \exp_{F_p(z)}(tDF_p(-v) + w')$ where $w' = o(t)$ as $t \to 0$, we have also

\begin{align*}
& \exp(bB(F_p(y), F_p(x), f(\xi))) = \exp\left(b\left(B(F_p(y), F_p(z), f(\xi)) - B(F_p(x), F_p(z), f(\xi))\right)\right) \\
& = \exp(b( 2< DF_p(v), \overrightarrow{F_p(z)f(\xi)}> t + < w - w', \overrightarrow{F_p(z)f(\xi)}> \\
& + \frac{1}{2}d^2 B^{f(\xi)}(t DF_p(-v)+w', tDF_p(-v)+w') - \frac{1}{2}d^2 B^{f(\xi)}(t DF_p(v)+w, tDF_p(v)+w) \\
& + o(t^2))) \\
& = \exp\left(b\left( 2< DF_p(v), \overrightarrow{F_p(z)f(\xi)}> t + < w - w', \overrightarrow{F_p(z)f(\xi)}> + o(t^2)\right)\right) \\
& = 1 + b \left( 2< DF_p(v), \overrightarrow{F_p(z)f(\xi)}> t + < w - w', \overrightarrow{F_p(z)f(\xi)}> \right) \\
& + 2b^2 < DF_p(v), \overrightarrow{F_p(z)f(\xi)}>^2 t^2 + o(t^2) \\
\end{align*}

\medskip

It follows that when computing the integral
\begin{align*}
& \int_{\partial X} \exp(bB(F_p(y), F_p(x), f(\xi))) d\mu^x_p(\xi) \\
& = \int_{\partial X} \exp(bB(F_p(y), F_p(x), f(\xi))) \frac{d\mu^x_p}{d\mu^z_p}(\xi) d\mu^z_p(\xi), \\
\end{align*}
after multiplying out the above expansions for $\exp(bB(F_p(y), F_p(x), f(\xi)))$ and $\frac{d\mu^x_p}{d\mu^z_p}(\xi)$, when we
integrate with respect to $\mu^z_p$, then, neglecting terms which are $o(t^2)$, the other terms involving $w,w'$ (which are $o(t)$)
vanish because the measure $f_* \mu^z_p$ is balanced
at $F_p(z)$, while the integrals of the terms $p(a(v,\xi) - \int_{\partial X} a(v, \xi) d\mu^z_p(\xi))t$ and
$(c(v,\xi) - \int_{\partial X} c(v,\xi) d\mu^z_p(\xi))t^2$ vanish since $\mu^z_p$ is a probability measure. Also
$$
\int_{\partial X} a(v, \xi) d\mu^z_p(\xi) = \int_{\partial X} < v, \overrightarrow{z\xi} > d\mu^z_p(\xi)
$$
because $f_* \mu^z_p$ is balanced at $F_p(z)$.
Thus we are finally left with

\begin{align*}
& \int_{\partial X} \exp(bB(F_p(y), F_p(x), f(\xi))) d\mu^x_p(\xi) \\
& = 1 + p^2 \left(\int_{\partial X} < v, \overrightarrow{z\xi} > d\mu^z_p(\xi)\right)^2 t^2 \\
  & + 2b \cdot p \cdot \left(\int_{\partial X} < v, \overrightarrow{z\xi} > < DF_p(v), \overrightarrow{F_p(z)f(\xi)} > - < DF_p(v), \overrightarrow{F_p(z)f(\xi)} >^2 d\mu^z_p(\xi)\right)t^2  \\
& + 2b^2 \left(\int_{\partial X} < DF_p(v), \overrightarrow{F_p(z)f(\xi)} >^2 d\mu^z_p(\xi)\right) t^2 + o(t^2) \\
\end{align*}

Now by Proposition \ref{implicit},
\begin{align*}
& p \cdot \left(\int_{\partial X} < v, \overrightarrow{z\xi} > < DF_p(v), \overrightarrow{F_p(z)f(\xi)} > - < DF_p(v), \overrightarrow{F_p(z)f(\xi)} >^2 d\mu^z_p(\xi)\right) \\
& = \int_{\partial X} d^2 B^{f(\xi)}(DF_p(v), DF_p(v)) d\mu^z_p(\xi) \\
& \geq 0 \\
\end{align*}

and hence

\begin{align*}
& \int_{\partial X} \exp(bB(F_p(y), F_p(x), f(\xi))) d\mu^x_p(\xi) \\
& \geq 1 + p^2 \left(\int_{\partial X} < v, \overrightarrow{z\xi} > d\mu^z_p(\xi)\right)^2 t^2 \\
& + 2b^2 \left(\int_{\partial X} < DF_p(v), \overrightarrow{F_p(z)f(\xi)} >^2 d\mu^z_p(\xi)\right) t^2 + o(t^2) \\
\end{align*}

as required. $\diamond$

\medskip

\begin{prop} \label{estimate2} We have
$$
\frac{p^2}{2b^2} \left(\int_{\partial X} < v, \overrightarrow{z\xi} > d\mu^z_p(\xi) \right)^2 \leq \int_{\partial X} ||DF_p(v)^{\bot f(\xi)}||^2 d\mu^z_p(\xi)
$$
\end{prop}

\medskip

\noindent{\bf Proof:} It follows from Propositions \ref{estimate1} and \ref{maininequality} that
\begin{align*}
& \cosh(bd(F_p(x), F_p(y))) \\
& \geq \int_{\partial X} \exp(bB(F_p(y), F_p(x), f(\xi))) d\mu^x_p(\xi) \\
& \geq 1 + p^2 \left(\int_{\partial X} < v, \overrightarrow{x\xi} > d\mu^z_p(\xi)\right)^2 t^2 \\
& + 2b^2 \left(\int_{\partial X} < DF_p(v), \overrightarrow{F_p(z)f(\xi)} >^2 d\mu^z_p(\xi)\right) t^2 + o(t^2) \\
\end{align*}
Now as $t \to 0$ we have $d(F_p(x), F_p(y)) = 2t||DF_p(v)|| + o(t)$, thus
$$
\cosh(bd(F_p(x), F_p(y))) = 1 + 2 t^2 b^2 ||DF_p(v)||^2 + o(t^2).
$$
It follows from the preceding inequality that
\begin{align*}
& 2 b^2 ||DF_p(v)||^2 \\
& \geq p^2 \left(\int_{\partial X} < v, \overrightarrow{x\xi} > d\mu^z_p(\xi)\right)^2  \\
& + 2b^2 \left(\int_{\partial X} < DF_p(v), \overrightarrow{F_p(z)f(\xi)} >^2 d\mu^z_p(\xi)\right). \\
\end{align*}
Using
$$
||DF_p(v)||^2 = \int_{\partial X} < DF_p(v), \overrightarrow{F_p(z)f(\xi)} >^2 + ||DF_p(v)^{\bot f(\xi)}||^2 d\mu^z_p(\xi)
$$
and the previous inequality, we obtain the inequality stated in the proposition. $\diamond$

\medskip

For each $z \in X$, we now let $\mu^z_{\infty}$ be a weak limit as $p \to \infty$ of the probability measures $\mu^z_p$ on
$\partial X$.

\medskip

\begin{prop} \label{xbalanced} For all $v \in T_z X$, we have
$$
\int_{\partial X} < v, \overrightarrow{z\xi} > d\mu^z_{\infty}(\xi) = 0,
$$
i.e. the probability measure $\mu^z_{\infty}$ on $\partial X$ is balanced at $z \in X$.
\end{prop}

\medskip

\noindent{\bf Proof:} From Propositions \ref{ineq1} and \ref{estimate2}, for any $v \in T^1_z X$, we have
\begin{align*}
&\frac{p^2}{2b^2} \left(\int_{\partial X} < v, \overrightarrow{z\xi} > d\mu^z_p(\xi) \right)^2 \\
& \leq \int_{\partial X} ||DF_p(v)^{\bot f(\xi)}||^2 d\mu^z_p(\xi) \\
& \leq p \int_{\partial X} < v, \overrightarrow{z\xi} >^2 d\mu^z_p(\xi) \\
\end{align*}
and hence
\begin{align*}
\left(\int_{\partial X} < v, \overrightarrow{z\xi} > d\mu^z_p(\xi) \right)^2 & \leq \frac{2b^2}{p} \int_{\partial X} < v, \overrightarrow{z\xi} >^2 d\mu^z_p(\xi) \\
                                                                    & \leq \frac{2b^2}{p}. \\
\end{align*}
Since $\mu^z_p \to \mu^z_{\infty}$ weakly as $p \to \infty$ along some sequence, passing to the limit above gives
$$
\int_{\partial X} < v, \overrightarrow{z\xi} > d\mu^z_{\infty}(\xi) = 0
$$
as required. $\diamond$

\medskip

\begin{lemma} \label{derivformula} For $x \in X, y \in Y, \xi \in \partial X$ we have
$$
\frac{df_*\rho_x}{d\rho_y}(f(\xi)) = \exp(B(y, \pi(\phi(\overrightarrow{x\xi})), f(\xi)))
$$
\end{lemma}

\medskip

\noindent{\bf Proof:} Let $z = \pi(\phi(\overrightarrow{x\xi})) \in Y$, then by definition of $\phi$ we
have $\frac{df_*\rho_x}{d\rho_z}(f(\xi)) = 1$, so by the Chain Rule
\begin{align*}
\frac{df_*\rho_x}{d\rho_y}(f(\xi)) & = \frac{df_*\rho_x}{d\rho_z}(f(\xi)) \cdot \frac{d\rho_z}{d\rho_y}(f(\xi)) \\
                                   & = 1 \cdot \exp(B(y, z, f(\xi))) \\
\end{align*}
$\diamond$

\medskip

\begin{prop} \label{supportmax} Let $K \subset \partial X$ be the set where the function
$$
\xi \in \partial X \mapsto \frac{df_*\rho_z}{d\rho_{F(z)}}(f(\xi))
$$
attains its maximum value. Then the support of the measure $\mu^z_{\infty}$
is contained in $K$.
\end{prop}

\medskip

\noindent{\bf Proof:} Let $M = \sup_{\xi \in \partial X} \frac{df_*\rho_z}{d\rho_{F(z)}}(f(\xi))$.
Given $\xi \in \partial X - K$, we can choose $\epsilon > 0$ and a neighbourhood $U$ of $\xi$ such that
$\frac{df_*\rho_z}{d\rho_{F(z)}}(f(\eta)) \leq M(1 - \epsilon)$ for all $\eta \in U$. Let $\psi$ be
a continuous function on $\partial X$ with support contained in $U$ such that $0 \leq \psi \leq 1$ on
$\partial X$ and such that $\psi = 1$ on a neighbourhood $N$ of $\xi$.

\medskip

For $\eta \in \partial X$,
$$
|B(F_p(z), \pi(\phi(\overrightarrow{z\eta})), f(\eta)) -  B(F(z), \pi(\phi(\overrightarrow{z\eta})), f(\eta))| \leq d(F_p(z), F(z)) \to 0
$$
as $p \to \infty$, so we can choose $p_0 > 1$ such that for $p \geq p_0$, by Lemma \ref{derivformula}
for all $\eta \in U$ we have
\begin{align*}
\exp(B(F_p(z), \pi(\phi(\overrightarrow{z\eta})), f(\eta))) & \leq \exp(B(F(z), \pi(\phi(\overrightarrow{z\eta})), f(\eta))) (1 + \epsilon) \\
                                                 & = \frac{df_*\rho_z}{d\rho_{F(z)}}(f(\eta)) (1 + \epsilon) \\
                                                 & \leq M(1 - \epsilon^2) \\
\end{align*}
while convergence of $L^p(\mu)$ norms to the $L^{\infty}(\mu)$ norm as $p \to \infty$ implies that
for $p_0$ chosen large enough, for $p \geq p_0$ we have
$$
c_{z,p} = ||\frac{df_*\rho_z}{d\rho_{F(z)}}(f(.))||^p_{L^p(\mu)} \geq (M(1 - \epsilon^2/2))^p
$$
hence
\begin{align*}
\int_{\partial X} \psi(\eta) d\mu^z_p(\eta) & \leq \int_{U} \exp(pB(F_p(z), \pi(\phi(\overrightarrow{z\eta})), f(\eta))) d\mu(\eta) \cdot c^{-1}_{z,p} \\
                                            & \leq \left(\frac{M(1 - \epsilon^2)}{M(1 - \epsilon^2/2)}\right)^p \\
                                            & \to 0 \\
\end{align*}
as $p \to \infty$. Letting $p \to \infty$ gives $\int_{\partial X} \psi(\eta) d\mu^z_{\infty}(\eta) = 0$, so
$\mu^z_{\infty}(N) = 0$ since $\psi = 1$ on $N$. It follows that $\xi$ is not in the support of $\mu^z_{\infty}$. $\diamond$

\medskip

\begin{prop} \label{inverses} Let $g : \partial Y \to \partial X$ be the inverse of the Moebius map $f : \partial X \to \partial Y$.
Then the circumcenter extensions $F : X \to Y$ and $G : Y \to X$ of the maps $f, g$ are inverses of each other.
\end{prop}

\medskip

\noindent{\bf Proof:} Let $x \in X$ and $y = F(x) \in Y$. Let $i_x : \partial X \to \partial X$
 denote the antipodal map of $\partial X$ centered
at the point $x$, i.e. the conjugate of the flip map $T^1_x X \to T^1_x X, v \mapsto -v$,
by the natural map $p : T^1_x X \to \partial X$.
Let $K \subset \partial X$ be the set where the function
$$
\xi \in \partial X \mapsto \log \frac{df_*\rho_x}{d\rho_y}(f(\xi)) = \log \frac{d\rho_x}{df^*\rho_y}(\xi)
$$
attains its maximum value. Then by Lemma \ref{maxmin}, $i_x(K)$ is contained in the set where the same function
attains its minimum value. By the Chain Rule, for $\xi \in \partial X$ and $\eta = f(\xi) \in \partial Y$, we have
$$
\log \frac{df_*\rho_x}{d\rho_y}(f(\xi)) + \log \frac{dg_*\rho_y}{d\rho_x}(g(\eta)) = 0.
$$
Letting $K' = f(i_x(K)) \subset \partial Y$, it follows that $K'$ is contained in the set $J$ say where the function
$$
\eta \in \partial Y \mapsto \log \frac{dg_*\rho_y}{d\rho_x}(g(\eta))
$$
attains its maximum value. Let $\nu_1$ be the probability measure $\nu_1 := (i_x)_* \mu^x_{\infty}$ on
$\partial X$ and let $\nu_2$ be the probability measure $\nu_2 = f_* \nu_1$ on $\partial Y$. By
Proposition \ref{supportmax} the support of $\mu^x_{\infty}$ is contained in $K$,
hence the support of $\nu_2$ is contained in $K'$. Now for any $\xi \in \partial X$,
we have $\overrightarrow{x i_x(\xi)} = - \overrightarrow{x\xi}$, and hence, for any $v \in T_x X$,
\begin{align*}
\int_{\partial X} < v, \overrightarrow{x\xi} > dg_* \nu_2(\xi) & = \int_{\partial X} < v, \overrightarrow{x\xi} > d\nu_1(\xi) \\
                                                    & = \int_{\partial X} < v, \overrightarrow{x i_x(\xi)} > d\mu^x_{\infty}(\xi) \\
                                                    & = \int_{\partial X} < v, - \overrightarrow{x\xi} > d\mu^x_{\infty}(\xi) \\
                                                    & = - \int_{\partial X} < v, \overrightarrow{x\xi} > d\mu^x_{\infty}(\xi) \\
                                                    & = 0 \\
\end{align*}
by Proposition \ref{xbalanced}. It follows that the measure $g_* \nu_2$ is balanced at $x \in X$, therefore,
since the support of $\nu_2$ is contained in the set $K' \subset J$, by (4) of
Proposition \ref{circbal} applied to the Moebius map $g : \partial Y \to \partial X$ we have $x = G(y)$.
$\diamond$

\medskip

\begin{prop} \label{circumlipschitz} The circumcenter extension $F : X \to Y$ is $\sqrt{b}$-Lipschitz.
\end{prop}

\medskip

\noindent{\bf Proof:} Let $x, y \in X$. For $\xi \in \partial X$, let
\begin{align*}
u_x(\xi) & = \log \frac{df_*\rho_x}{d\rho_{F(x)}}(f(\xi)) = \exp(B(F(x), \overrightarrow{x\xi}, f(\xi))) \\
u_y(\xi) & = \log \frac{df_*\rho_y}{d\rho_{F(y)}}(f(\xi)) = \exp(B(F(y), \overrightarrow{y\xi}, f(\xi))) \\
\end{align*}
and let $K_x, K_y$ be the subsets of $\partial X$ where the functions $u_x, u_y$ attain their maximum values
respectively, and let the maximum values be $r(x), r(y)$ respectively. Then for $\xi \in K_x, \eta \in K_y$,
we have
\begin{align*}
& B(F(y), F(x), f(\xi)) \\
& = B(F(y), \pi(\phi(\overrightarrow{y\xi})), f(\xi)) + B(\pi(\phi(\overrightarrow{y\xi})), \pi(\phi(\overrightarrow{x\xi})), f(\xi)) +
B(\pi(\phi(\overrightarrow{x\xi})), F(x), f(\xi)) \\
& = u_y(\xi) + B(y, x, \xi) - r(x) \\
\end{align*}
and similarly
$$
B(F(x), F(y), f(\eta)) = u_x(\eta) + B(x, y, \eta) - r(y),
$$
thus
\begin{align*}
& B(F(y), F(x), f(\xi)) + B(F(x), F(y), f(\eta)) \\
& = (B(y, x, \xi) + B(x, y, \eta)) + (u_y(\xi) - r(y)) + (u_x(\eta) - r(x)) \\
& \leq B(y, x, \xi) + B(x, y, \eta). \\
\end{align*}
Since the measures $f_* \mu^x_p, f_* \mu^y_p$ are balanced at the points $F_p(x), F_p(y)$ respectively,
and $F_p(x) \to F(x), F_p(y) \to F(y)$ as $p \to \infty$, it is easy to see that the measures
$f_* \mu^x_{\infty}, f_* \mu^y_{\infty}$ are balanced
at the points $F(x), F(y)$ respectively, so an argument similar to the proof of Proposition \ref{maininequality}
(using Lemma \ref{buseform}) gives
\begin{align*}
\cosh(d(F(x), F(y))) \leq \int_{\partial X} \exp(B(F(y), F(x), f(\xi))) d\mu^x_{\infty}(\xi) \\
\cosh(d(F(y), F(x))) \leq \int_{\partial X} \exp(B(F(x), F(y), f(\eta))) d\mu^y_{\infty}(\eta). \\
\end{align*}
The measures $\mu^x_{\infty}, \mu^y_{\infty}$ are also balanced at the points $x, y$ respectively, so
again an argument similar to the proof of Proposition \ref{maininequality} gives
\begin{align*}
\cosh(bd(x, y)) \geq \int_{\partial X} \exp(bB(y, x, \xi)) d\mu^x_{\infty}(\xi) \\
\cosh(bd(y, x)) \geq \int_{\partial X} \exp(bB(x, y, \eta)) d\mu^y_{\infty}(\eta). \\
\end{align*}
Now using the fact that the supports of $\mu^x_{\infty}, \mu^y_{\infty}$ are contained in $K_x, K_y$
respectively, we have
\begin{align*}
& \cosh^2(d(F(x), F(y))) \\
& \leq \left(\int_{\partial X} \exp(B(F(y), F(x), f(\xi))) d\mu^x_{\infty}(\xi)\right)
\left(\int_{\partial X} \exp(B(F(x), F(y), f(\eta))) d\mu^y_{\infty}(\eta)\right) \\
& = \int_{\partial X} \int_{\partial X} \exp(B(F(y), F(x), f(\xi)) + B(F(x), F(y), f(\eta))) d\mu^x_{\infty}(\xi) d\mu^y_{\infty}(\eta) \\
& \leq \int_{\partial X} \int_{\partial X} \exp(B(y, x, \xi) + B(x, y, \eta)) d\mu^x_{\infty}(\xi) d\mu^y_{\infty}(\eta) \\
& = \left(\int_{\partial X} \exp(B(y, x, \xi)) d\mu^x_{\infty}(\xi)\right)
\left(\int_{\partial X} \exp(B(x, y, \eta)) d\mu^y_{\infty}(\eta)\right) \\
& \leq \left(\int_{\partial X} \exp(bB(y, x, \xi)) d\mu^x_{\infty}(\xi)\right)^{1/b}
\left(\int_{\partial X} \exp(bB(x, y, \eta)) d\mu^y_{\infty}(\eta)\right)^{1/b} \\
& \leq (\cosh^2(bd(x,y)))^{1/b} \\
\end{align*}
thus
$$
\cosh^{b}(d(F(x), F(y))) \leq \cosh(bd(x,y)).
$$
Now for $t \geq 0$, $\cosh^b(t) \geq 1 + bt^2/2$, and there is a
universal constant $C > 0$ such that $\cosh bt \leq 1 + b^2 t^2$ for $0 \leq t \leq C$,
so $d(x,y) \leq C$ implies
$$
\frac{1}{2} d(F(x), F(y))^2 \leq b d(x, y)^2,
$$
thus $F$ is locally Lipschitz. It follows that $F$ is differentiable almost everywhere. At a point $x$
of differentiability of $F$, for $v \in T^1_x X$ letting $y = \exp_x(tv)$, a Taylor expansion of both sides
of the inequality $\cosh^{b}(d(F(x), F(y))) \leq \cosh(bd(x,y))$ up to second order in $t$ easily gives
$||DF(v)|| \leq \sqrt{b}$. Now it is a standard fact that if $F$ is a locally Lipschitz map between
complete Riemannian manifolds such that $||DF|| \leq C$ almost everywhere then $F$ is $C$-Lipschitz, hence
in our case $F$ is $\sqrt{b}$-Lipschitz. $\diamond$

\medskip

\noindent{\bf Proof of Theorem \ref{mainthm}:} Since the maps $F : X \to Y, G : Y \to X$ are inverses
of each other, it follows from Proposition \ref{circumlipschitz} above that they are $\sqrt{b}$-bi-Lipschitz
homeomorphisms.
$\diamond$

\medskip

\noindent{\bf Proof of Theorem \ref{rigid}:} Given $X, Y$ closed, negatively curved $n$-manifolds, after 
rescaling the metrics on $X, Y$ by the same constant $C > 0$, we may assume by choosing $C$ appropriately that 
both manifolds have sectional curvatures bounded above by $-1$. Then the universal covers $\tilde{X}, \tilde{Y}$ 
are CAT(-1) spaces and equality of the marked length spectra of $X, Y$ implies existence of an equivariant 
Moebius homeomorphism $f : \partial\tilde{X} \to \partial\tilde{Y}$. By the naturality of the circumcenter extension, 
the circumcenter extension $\hat{f} : \tilde{X} \to \tilde{Y}$ is equivariant, and is a bi-Lipschitz homeomorphism by 
the previous theorem, hence induces a bi-Lipschitz homeomorphism $F : X \to Y$. $\diamond$ 
\bibliography{moeb}
\bibliographystyle{alpha}

\end{document}